\UseRawInputEncoding
\documentclass[10pt]{article}
\usepackage{amsmath}
\numberwithin{equation}{section}
\usepackage{amsfonts,amssymb,amsmath}
\usepackage{mathrsfs}
\usepackage{color}
\newcommand{\beq}{\begin{equation}}
\newcommand{\enq}{\end{equation}}

\newtheorem{Theorem}{Theorem}[section]
\newtheorem{Proposition}{Proposition}[section]
\newtheorem{Lemma}[Theorem]{Lemma}

\newtheorem{Remark}[Theorem]{Remark}

\newcommand{\benu}{\begin{enumerate}}
\newcommand{\beqa}{\begin{eqnarray}}
\newcommand{\beqan}{\begin{eqnarray*}}
\newcommand{\eay}{\end{array}}
\newcommand{\edm}{\end{displaymath}}
\newcommand{\eenu}{\end{enumerate}}
\newcommand{\eeq}{\end{equation}}
\newcommand{\eeqa}{\end{eqnarray}}
\newcommand{\eeqan}{\end{eqnarray*}}

\newcommand{\br}{\begin{Remark}}
\newcommand{\er}{\end{Remark}}

\newcommand{\bqa}{\begin{eqnarray}}
\newcommand{\eqa}{\end{eqnarray}}
\newcommand{\bqw}{\begin{eqnarray*}}
\newcommand{\eqw}{\end{eqnarray*}}

\newcommand{\bea}{\begin{array}{cc}}
\newcommand{\ena}{\end{array}}

\allowdisplaybreaks[4]
\begin{document}
\begin{center}
{\large \bf  Local existence of solutions to 3D Prandtl equations with a special structure}\\
\vspace{0.25in} ${\rm \rm Yuming\; Qin^{1*}}\quad{\rm Xiuqing\; Wang^{2  }}   $\\
\vspace{0.15in} 1*. Department of Mathematics, Institute for Nonlinear Sciences,
Donghua University,\\
Shanghai 201620, P. R. China.\\
E-mail: yuming\_qin@hotmail.com,\quad yuming@dhu.edu.cn\\

2. Department of Mathematics, Kunming University of Science and Technology, \\ Kunming 650500,
Yunnan, P. R. China.\\
E-mails:  daqingwang@kust.edu.cn

 \vspace{3mm}
\end{center}
\vspace{0.2in}

\begin{abstract}~~
In this paper, we consider the 3D   Prandtl equation in a periodic domain and prove the local existence and uniqueness of solutions by the energy method in a polynomial weighted  Sobolev space. Compared to the existence and uniqueness of solutions to the classical Prandtl equations
where the Crocco transform has always been used with the general outer flow $U\neq\text{constant}$, this
Crocco transform is not needed here for 3D  Prandtl equations.
 We use the skill of cancellation mechanism and construct a new unknown function to show that   the
 existence and uniqueness of solutions to 3D Prandtl equations (cf.  Masmoudi and Wong, Comm. Pure Appl. Math.,
    68(10)(2015), 1683-1741) which extends from the two dimensional case in  \cite{12}
to the present three dimensional case  with a
special structure.


 \end{abstract}

\section{Introduction}
Experimental data and theoretical analysis show that in many important practical cases, for fluids whose viscosity is small, there exists a thin transition layer near the boundary, in which the behavior of flow changes dramatically, this phenomenon is called boundary layer theory in \cite{2}.
In this paper, we consider the 3D Prandtl equations in the periodic domain  $\{(t,x,y,z)\big|t>0, (x,y)\in \mathbb{T}^{2}, z\in\mathbb{R^{+}}\}$:
\begin{equation}\left\{
\begin{array}{ll}
\partial _{t}u+(u\partial _{x} +v\partial _{y} +w\partial_{z})u+\partial _{x}P=\partial _{z}^{2}u,\\
\partial _{t}v+(u\partial _{x} +v\partial _{y} +w\partial_{z})v+\partial _{y}P=\partial _{z}^{2}v,\\
\partial _{x}u+\partial _{y}v+\partial_{z}w=0,\\
(u,v)|_{t=0}=(u_{0}(x,y,z),v_{0}(x,y,z)),\\
(u,v,w)|_{z=0}=0,\\
  \lim\limits_{z\rightarrow+\infty}(u,v)=(U(t,x,y), V(t,x,y)).
\end{array}
 \label{1.1}         \right.\end{equation}
Here $\mathbb{T}^{2}\subseteq \mathbb{R}^{2}$,
$(u,v,w)=(u(t,x,y,z),v(t,x,y,z),w(t,x,y,z))$ denote the velocity field, $(U(t,x,y), V(t,x,y))$ and $P$ are the traces of the tangential velocity field and the pressure of the the Euler flow respectively which satisfy Bernoulli��s law
\begin{equation}\left\{
\begin{array}{ll}
\partial _{t}U+(U\partial _{x} +V\partial _{y} )U+\partial _{x}P=0,\\
\partial _{t}V+(U\partial _{x} +V\partial _{y}  )V+\partial _{y}P=0.
\end{array}
 \label{1.2}         \right.\end{equation}

Despite of its importance in physics, there are very few mathematical results
on the Prandtl equations in three space variables, and the three-dimensional results are based on special structures, such as \cite{4,5}, where authors constructed a solution of the three-dimensional Prandtl equations $(\ref{1.1})$ with the structure
\begin{eqnarray}
 (u(t,x,y,z),K(t, x,y)u(t,x,y,z),w(t,x,y,z)),
 \label{1.3}
\end{eqnarray}
and the outer Euler flow takes the following form on the boundary $\{z=0\}$,
$$(U(t,x,y), K(t, x,y)U(t,x,y),0 ). $$
But in what follows, we shall consider the equivalent system of system (\ref{1.1}) (see Proposition \ref{p1.2}  for the specific proof). In fact, the regularized equivalent system of 3D Prandtl equations (\ref{1.1}) can be expressed as follows,
\begin{equation}\left\{
\begin{array}{ll}
\partial _{t}u^\epsilon+(u^\epsilon\partial _{x} +K^\epsilon u^\epsilon\partial _{y} +w^\epsilon\partial_{z})u^\epsilon -\epsilon^2\partial _{x}^{2}u^\epsilon-\epsilon^2\partial _{y}^{2}u^\epsilon=\partial _{z}^{2}u^\epsilon-\partial_x P^\epsilon,\\
\partial _{x}u^\epsilon+\partial _{y}(K^\epsilon u^\varepsilon)+\partial_{z}w^\epsilon=0,\\
 u ^\epsilon \big|_{t=0}= u_{0}(x,y,z) ,\\
(u^\epsilon,K^{^\epsilon}u^\epsilon,w^\epsilon)\big|_{z=0}=0,\quad  \lim\limits_{z\rightarrow+\infty}u^\epsilon=U(t, x, y)
\end{array}
 \label{1.4}         \right.\end{equation}
with  regularized Bernoulli's law
\begin{eqnarray}
\partial_t U+ U\partial_x U+ K^\epsilon U\partial_y U-\epsilon^2 \partial_x^2 U -\epsilon^2 \partial_y^2 U +\partial_x P^\epsilon=0.
\label{1.44}
\end{eqnarray}
Let the vorticity $\varphi^\epsilon=\partial_{z}u^\epsilon$, then equations $(\ref{1.4})$ satisfy the following vorticity system: for any $\epsilon >0$,
\begin{equation}\left\{
\begin{array}{ll}
\partial _{t}\varphi^\epsilon+(u^\epsilon\partial _{x} +K^\epsilon u^\epsilon\partial _{y}+w^\epsilon\partial_{z})\varphi^\epsilon -\varepsilon\partial _{x}^{2}\varphi^\epsilon-\varepsilon\partial _{y}^{2}\varphi^\epsilon= \frac{1}{2}\partial_y K \partial_z |u^\epsilon|^2 +\partial _{z}^{2}\varphi^\epsilon,\\
\partial _{x}u^\epsilon+\partial _{y}(K^\epsilon u^\epsilon)+\partial_{z}w^\epsilon=0,\\
\varphi^\epsilon(0,x,y,z)=\partial_{z}u_{0},\\
\partial _{z}\varphi^\epsilon|_{z=0}=\partial_x P^\epsilon,
\end{array}
 \label{1.5}         \right.\end{equation}
where the velocity field $(u^\epsilon, w^\epsilon)$ is given by
\begin{eqnarray}
u^\epsilon (t, x, y, z) = U-\int_{z}^{+\infty} \varphi^\epsilon (t, x, y, \tilde{z})d\tilde{z},
\end{eqnarray}
and
\begin{eqnarray}
w^\epsilon (t, x, y, z) = -\int_{0}^{z} \partial_x u^\epsilon (t, x, y, \tilde{z})d\tilde{z}-\int_{0}^{z} \partial_y (Ku^\epsilon (t, x, y, \tilde{z}))d\tilde{z}.
\end{eqnarray}
Next, we introduce the weighted Sobolev space and define the space $H^{s,\gamma}_{\sigma,\delta}$ by
$$H^{s,\gamma}_{\sigma,\delta}:=\left\{\varphi:\mathbb{T}^{2}\times\mathbb{R^{+}}\rightarrow\mathbb{R},
\|\varphi\|_{H^{s,\gamma}}< +\infty, (1+z)^{\sigma}|\varphi| \geq\delta,
\sum \limits_{|\alpha|\leq 2}  |(1+z)^{\sigma+\alpha_{3}}D^{\alpha}\varphi | ^{2}\leq\frac{1}{\delta^{2}} \right\},$$
where $D^{\alpha}:=\partial_{x}^{\alpha_{1}}\partial_{y}^{\alpha_{2}}\partial_{z}^{\alpha_{3}}$, $\alpha_{1}+\alpha_{2}+\alpha_{3}=s$,  $s\geq 5$, $\gamma\geq 1$, $\sigma>\gamma+\frac{1}{2}$ and $\delta\in (0,1)$. We define the norm as
 $$\|\varphi\|_{H^{s,\gamma}}^{2}:=\sum \limits_{|\alpha|\leq s} \|(1+z)^{\gamma+\alpha_{3}}D^{\alpha}\varphi\|_{L^{2}}^{2}$$
and
$$\|\varphi\|_{H^{s,\gamma}_{g}}^{2}:=\|(1+z)^{\gamma}g_{s}\|^{2}_{L^{2}}+
\sum\limits_{\substack{ |\alpha|\leq s\\ \alpha_{1}+\alpha_{2}\leq s-1}}  \|(1+z)^{\gamma+\alpha_{3}}D^{\alpha}\varphi\|_{L^{2}}^{2}.$$
Here
$$g_{s}:=\partial_{xy}^{s}\varphi-\frac{\partial_{z}\varphi}{\varphi}\partial_{xy}^{s}(u-U)
=\sum\limits_{i=0}^{ s}\left( \partial_{x}^{i}\partial_{y}^{s-i}\varphi- \frac{\partial_{z}\varphi}{\varphi}  \partial_{x}^{i}\partial_{y}^{s-i} (u-U) \right),$$
and
$$g_{s}=(g_{s})_{x}+(g_{s})_{y}=\partial_x\partial_{xy}^{s-1}\varphi-\frac{\partial_{z}\varphi}{\varphi}\partial_x\partial_{xy}^{s-1}(u-U)+\partial_y\partial_{xy}^{s-1}\varphi
-\frac{\partial_{z}\varphi}{\varphi}\partial_y\partial_{xy}^{s-1}(u-U),$$
if $\varphi >0$.
Hence, we can calculate that
$$g_{j+1}=\partial_{x}g_{j}+\partial_{y}g_{j}+\partial_{xy}^{j}u\partial_{x}a+\partial_{xy}^{j}u\partial_{y}a,\quad
a=\frac{\partial_{z}\varphi}{\varphi}.$$

In this paper, for convenience, we simply write
$$\partial_{xy}^{s}:= \sum\limits_{i=0}^{ s} \partial_{x}^{i}\partial_{y}^{s-i} , \quad  \iiint\cdot:=\iint_{\mathbb{T}^{2}}\int_{\mathbb{R^{+}}}\cdot dxdydz.$$

Next, we state our main result as follows.
\begin{Theorem}\label{t1.1}
Given any even integer  $s \geq 5$, and real numbers $\gamma, \sigma, \delta$ satisfying $\gamma\geq 1$, $\sigma>\gamma+\frac{1}{2}$ and $\delta\in (0,1)$.   Assume the following conditions on the initial data,  the outer flow $U$ and $K(x,y)$:

(i)~Suppose that the initial data $u_0-U(0,x) \in H^{s, \gamma-1}$ and $\partial_z u_0 \in H^{s,\gamma}_{\sigma, 2\delta}$ satisfy the compatibility conditions $u_0|_{z=0}$ and $\lim \limits_{z \rightarrow +\infty} u=U|_{t=0}$. In addition, when $s = 5$, it is further assumed that $\delta \geq 0$ is chosen small enough such that $\|\omega_0\|_{H^{s,\gamma} } \leq C \delta^{-1}$ with a generic constant $C$.

(ii)~The outer flow $U$ is supposed to satisfy
 \begin{eqnarray}
\sup \limits_{t} \sum\limits_{l=0}^{\frac{s}{2}+1}\|\partial_t^l U\|_{H^{s-2l+2}(\mathbb{T}^2)} < + \infty.
\label{1.8}
\end{eqnarray}

(iii)~The $K(x,y)$ is supposed to satisfy
 \begin{eqnarray}
\|\partial^{s}_{xy} K\|_{L^{\infty}(\mathbb{T}^2)} < + \infty.
\label{1.88}
\end{eqnarray}

Then there exist a time $T := T(K, s, \gamma, \sigma, \delta,  \|\varphi_0\|_{H^{s,\gamma} },U)$ such that the initial-boundary value problem
(\ref{1.1})-(\ref{1.3}) has a unique classical solution $(u, v, w)$ (v=Ku) satisfying
$$u-U \in L^{\infty}([0,T];H^{s,\gamma-1}) \cap C([0,T]; H^{s}-\varphi)$$
and
$$ \partial_z u \in L^{\infty}([0,T];H^{s,\gamma}_{\sigma, \delta}) \cap C([0,T]; H^{s}-\varphi),$$
where $H^s-\varphi$ is the space $H^s$ endowed with its weak topology.
\end{Theorem}

Let us now briefly review the background and corresponding results about the boundary layer. The mathematical studies on the Prandtl boundary layer have a very long history. Moore \cite{8} gave an analytical framework of the three-dimensional boundary layer in 1956.
In the following decades, the 2-dimensional boundary layer theory developed rapidly.
The first well-known result was developed by Oleinik and Samokhin in \cite{2}, where under the
monotonicity condition on tangential velocity with respect to the normal variable to the boundary, the local (in time) well-posedness of Prandtl equations was obtained by using the Crocco transformation and von Mises transformation under the outer flow $U\neq\text{constant}$.  Since then, Crocco transformation and von Mises transformation have been used in boundary layer problems. For example, Xin and Zhang \cite{7} established a global existence of weak solutions to the two-dimensional Prandtl system for the pressure is favourable (i.e., $\partial_{x}P<0$), which had generalized the local well-posedness results of Oleinik \cite{2}.
But, Masmoudi and Wong  \cite{3} proved the local existence and uniqueness of solutions for  classical Prandtl  system by energy methods in a polynomial weighted  Sobolev space under the outer flow $U\neq\text{constant}$,  and Fan, Ruan, Yang \cite{11} proved the  local well-posedness for the compressible Prandtl boundary layer equations by the same method.
The authors Qin, Wang and Liu \cite{12} proved the local existence and uniqueness of the Prandtl-Hartmann regime with  the general outer flow $U\neq\text{constant}$, and Dong and Qin \cite{13} established the global well-posedness of solutions in $H^s(1\leq s\leq 4)$ with $ U= \text{constant}$ and without monotonicity condition and a lower bound.

Compared to the 2D case, the results of the three-dimensional boundary layer equations were very few.
Until recent years, Liu, Yang and Wang \cite{4,5} obtained the local and global existence of weak solutions to the three-dimensional Prandtl equations with a special structure by Crocco transformation.
Liu, Yang and Wang \cite{14} gave an ill-posedness criterion for the Prandtl equations in three-dimensional space  under the outer flow  $(U,V)=\text{constant}$, which shows that the monotonicity condition on tangential velocity fields is not sufficient for the well-posedness of the three-dimensional Prandtl equations. Based on this result,
they proved the existence and uniqueness of local solutions of three-dimensional boundary layer equations without any special structure in the Gevrey function space, see \cite{15,17}.
Lin and Zhang \cite{16}  proved the almost global existence of classical solutions to the
3D Prandtl system with the initial data which lie within $\varepsilon$ of a stable shear flow under the outer flow  $(U,V)=\text{constant}$.
Pan and Xu \cite{18} proved the global existence of solutions to the three dimensional axially
symmetric Prandtl boundary layer equations with small initial data under $(U,V)=0$.

As we have known,   Crocco transformation  and the outer flow  $(U,V)=\text{constant}$  have been used in three-dimensional prandtl system, such as \cite{15,17,16,14,4,5,18,19}.
But, we need neither the Crocco transformation nor the the outer flow  $(U,V)=\text{constant}$ in this paper. In fact, we  first use the energy method to prove the local existence of three-dimensional Prandtl system based on a cancellation property. Our result is an improvement of the   existence and  uniqueness of solutions   in Masmoudi and Wong \cite{3}  from the two dimensional case to the three dimensional with special structure case.

The paper is arranged as follows. Section 1 is an introduction, which contains many inequalities used in this paper. In Section 2, we shall give the uniform estimates with the weighted norm. The difficulty of this paper is the estimates on the boundary at $z=0$, and we give the regular pattern at $z=0$ for $s\geq 5$ eventually through a series of precise calculations, see Lemmas \ref{y2.1}-\ref{y2.2}. In Section 3, we shall prove local existence and uniqueness of solutions to the Prandtl system to problem (\ref{1.1}).
\subsection{Preliminaries}
\begin{Proposition} \label{p1.2}
Under the assumption $\partial_z u >0$ and special structures $v(t,x,y,z)=K(t, x, y)u(t,x,y,z)$, to study
the problem (\ref{1.1}) is equivalent to studying the following reduced problem:
\begin{equation}\left\{
\begin{array}{ll}
\partial _{t}u+(u\partial _{x} +K u\partial _{y} +w\partial_{z})u =\partial _{z}^{2}u^\epsilon-\partial_x P,\\
\partial _{x}u+\partial _{y}(K u)+\partial_{z}w=0,\\
 u  \big|_{t=0}= u_{0}(x,y,z) ,\\
(u,w)\big|_{z=0}=0,\quad  \lim\limits_{z\rightarrow+\infty}u=U(t, x, y),
\end{array}
 \label{4.40}         \right.\end{equation}
 where  $K$ is independent of $t$, and $U$ satisfies
 \begin{equation}\left\{
\begin{array}{ll}
\partial_t U+ U\partial_x U+ K U\partial_y U+\partial_x P=0,\\
\partial_y P-K\partial_x P=0.
\end{array}
        \right.\end{equation}
        \end{Proposition}
$\mathbf{Proof.}$
Using special structures $v(k,x,y,z)=K(t, x, y)u(k,x,y,z)$ for the equation $(\ref{1.1})_2$, we have
 \begin{eqnarray}
\partial_t (Ku)+(u \partial_x +Ku \partial_y +w \partial_z)(Ku)+\partial_y P - K\partial_z^2 u=0,
\label{4.41}
\end{eqnarray}
which, by the equation $(\ref{1.1})_1$, yields
 \begin{eqnarray}
u\left[\partial_t K+u( \partial_x + \partial_y)K\right]-K\partial_x P + \partial_y P=0.
\label{4.42}
\end{eqnarray}
Noting that $K(t,x,y)$ is independent of $z$, by differentiating (\ref{4.42}) with respect to $z$, it follows that
 \begin{eqnarray}
\partial_z u\partial_t K+2 u \partial_z u (\partial_x+K\partial_y)K=0.
\label{4.43}
\end{eqnarray}
Due to the monotonicity assumption, we can further get
 \begin{eqnarray}
\partial_t K+2 u  (\partial_x+K\partial_y)K=0.
\label{4.44}
\end{eqnarray}
Differentiating (\ref{4.44}) with respect to $z$ gives the $Burgers$ equation
 \begin{eqnarray}
  (\partial_x+K\partial_y)K=0.
\label{4.45}
\end{eqnarray}
Plugging (\ref{4.45}) into (\ref{4.44}), it follows
 \begin{eqnarray}
 \partial_t K=0,
\label{4.46}
\end{eqnarray}
which, combined (\ref{4.42}) with (\ref{4.45}), gives
 \begin{eqnarray}
-K\partial_x P + \partial_y P=0.
\label{4.47}
\end{eqnarray}
The above equation means that $(\partial_x P, \partial_y P)$ is parallel to both the velocity field of the outer Euler flow and the tangential velocity field in the boundary layer. Now, we assume that  the classical solution $(u, v, w)$  to the problem (\ref{1.1}) satisfies $v(t,x,y,z)=K(x,y)u(t,x,y,z)$, then  by using (\ref{4.45}) and (\ref{4.47}), $W(t, x, y, z):\equiv v(t,x,y,z)-K(x,y)u(t,x,y,z)$ satisfies the following problem:
\begin{equation}\left\{
\begin{array}{ll}
\partial _{t}W+(u\partial _{x} +v\partial _{y} +w\partial_{z})W+(\partial_y v - K\partial_y u)W=\partial _{z}^{2}W,\\
W|_{t=0}=0,\\
W|_{z=0}=0,\\
  \lim\limits_{z\rightarrow+\infty}W=0.
\end{array}
 \label{4.48}         \right.\end{equation}
The equation (\ref{4.48}) has only trivial solution $W\equiv0$ based on the energy argument. Therefore,   studying
the problem (\ref{1.1}) is equivalent to studying  the following reduced problem for only two
unknown functions $u$ and $w$,
\begin{equation}\left\{
\begin{array}{ll}
\partial _{t}u+(u\partial _{x} +K u\partial _{y} +w\partial_{z})u =\partial _{z}^{2}u^\epsilon-\partial_x P,\\
\partial _{x}u+\partial _{y}(K u)+\partial_{z}w=0,\\
 u  \big|_{t=0}= u_{0}(x,y,z) ,\\
(u,w)\big|_{z=0}=0,\quad  \lim\limits_{z\rightarrow+\infty}u=U(t, x, y) .
\end{array}
 \label{4.49}         \right.\end{equation}
Hence, the proof of   Proposition \ref{p1.2} is complete.
$\hfill$$\Box$

\begin{Proposition}
Let $s\geq 5$ be an integer, $\gamma\geq 1$, $\sigma>\gamma+\frac{1}{2}$ and $\delta\in (0,1)$. Then for any $\varphi \in H^{s,\gamma}_{\sigma,\delta}(\mathbb{T}^{2}\times\mathbb{R}^{+})$,
we have the following inequality
 \begin{eqnarray}
C_{\delta }\|\varphi\|_{H^{s,\gamma}_{g}}\leq  \|\varphi\|_{H^{s,\gamma} }+ \|u-U\|_{H^{s,\gamma-1} }
\leq C_{s, \gamma,\sigma,\delta  }\left(\|\varphi\|_{H^{s,\gamma}_{g}}+\|\partial_{xy}^s U\|_{L^2(\mathbb{T}^2)}\right)
  ,
\label{4.1}
\end{eqnarray}
where $C_{s, \gamma,\sigma,\delta  }>0$ is a constant and only depends on $s, \gamma,\sigma,\delta  $.
\end{Proposition}

\begin{Lemma} \label{y4.1}
Let $f:\mathbb{T}^{2} \times\mathbb{R}^{+}\rightarrow \mathbb{R}$, \\
 (i)  if $\lambda>-\frac{1}{2}$ and $\lim\limits_{z\rightarrow+\infty}f(x,y ,z)=0$, then
  \begin{eqnarray}
  \|(1+z)^{\lambda}f\|_{L^{2}(\mathbb{T}^{2} \times\mathbb{R}_{+})} \leq \frac{2}{2\lambda+1} \|(1+z)^{\lambda+1}\partial_{z}f\|_{L^{2}(\mathbb{T}^{2} \times\mathbb{R}_{+})};
 \label{4.10}
 \end{eqnarray}
 (ii)  if $\lambda<-\frac{1}{2}$ and $ f(x,y,z )|_{z=0}=0$, then
 \begin{eqnarray}
 \|(1+z)^{\lambda}f\|_{L^{2}(\mathbb{T}^{2} \times\mathbb{R}^{+})} \leq
 -\frac{2}{2\lambda+1} \|(1+z)^{\lambda+1}\partial_{z}f\|_{L^{2}(\mathbb{T}^{2} \times\mathbb{R}^{+})}.
\label{4.11}
\end{eqnarray}
\end{Lemma}
$\mathbf{Proof.}$
The proof is elementary.
Using partial integrating  in $z$-variable proves $(\ref{4.10})$-$(\ref{4.11})$.
$\hfill$$\Box$

\begin{Lemma} \label{y4.2}
Let $s\geq 5$ be an integer, $\gamma\geq 1$, $\sigma>\gamma+\frac{1}{2}$ and $\delta\in (0,1)$. If $\varphi \in H^{s,\gamma}_{\sigma,\delta}(\mathbb{T}^{2}\times\mathbb{R}^{+})$, then  for any $k=0,1,2,\cdots, s$,
 \begin{eqnarray}
 \|(1+z)^{\gamma} g_{k}\|_{L^{2}}\leq\|(1+z)^{\gamma}\partial _{xy}^{k}\varphi\|_{L^{2}}+\delta^{-2}\|(1+z)^{\gamma-1}\partial _{xy}^{k}(u-U)\|_{L^{2}}  ,
 \label{4.5}
\end{eqnarray}
where $g_{k}:=\partial_{xy}^{k}\varphi-\frac{\partial_{z}\varphi}{\varphi}\partial_{xy}^{k}(u-U)$. In addition, if $u \big|_{z=0} = 0$, then for any $k=1, 2, \ldots, s$,
 \begin{eqnarray}
 \|(1+z)^{\gamma}\partial _{xy}^{k}\varphi\|_{L^{2}}+\|(1+z)^{\gamma-1}\partial _{xy}^{k}(u-U)\|_{L^{2}}
\leq C_{s, \gamma,\sigma,\delta  }\left(\|(1+z)^{\gamma} g_k\|_{L^2}+\|\partial_{xy}^k U\|_{L^2(\mathbb{T}^2)}\right) , \label{4.4}
\end{eqnarray}
where $C_{s, \gamma,\sigma,\delta  }>0$ is a constant and only depends on $s, \gamma,\sigma,\delta  $.
\end{Lemma}
The proof is similar to that given in \cite{3}, it is only related to the defined weighted Sobolev space and norms, not to the condition $(\ref{1.4})_2$, so we omit it here for brevity.

\begin{Lemma}(\cite{20}) \label{y4.3}
Let $f:\mathbb{T}^{2} \times\mathbb{R}^{+}\rightarrow \mathbb{R}$ and $ f(x,y,z )|_{z=0}=0$. Then there exists a constant $C>0$ such that
\begin{eqnarray}
\begin{aligned}
&  \|f(x,y,z )\|_{L^{\infty}(\mathbb{T}^{2} \times \mathbb{R}^{+})}  \\
&   \leq C
 (\|f(x,y,z)\|_{L^{2}(\mathbb{T}^{2}\times \mathbb{R}^{+} )}   +\|\partial_{z} f(x,y,z)\|_{L^{2}(\mathbb{T}^{2}\times \mathbb{R}^{+} )}
  +\partial_{xy}^{1}\|f(x,y,z)\|_{L^{2}(\mathbb{T}^{2}\times \mathbb{R}^{+} )}\\
 &\quad +\|\partial_{z}\partial_{xy}^{1}f(x,y,z)\|_{L^{2}(\mathbb{T}^{2}\times \mathbb{R}^{+} )}
  +\|\partial_{xy}^{2}f(x,y,z)\|_{L^{2}(\mathbb{T}^{2}\times \mathbb{R}^{+} )}+\|\partial_{z}\partial_{xy}^{2}f(x,y,z)\|_{L^{2}(\mathbb{T}^{2}\times \mathbb{R}^{+} )}  ).
\label{4.12}
\end{aligned}
\end{eqnarray}

\end{Lemma}

$\mathbf{Proof.}$
Firstly, we have
\begin{eqnarray}
&& f^{2}(x,y,z)=\int_{0}^{z}\partial_{\tau}[f^{2}(x,y,\tau)]d\tau =2\int_{0}^{z}f (x,y,\tau)\partial_{\tau} f (x,y,\tau) d\tau   \nonumber \\
&& \leq C\left(\int_{0}^{z} f^{2}(x,y,\tau) d\tau \right)^{\frac{1}{2}} \left(\int_{0}^{z}| \partial_{\tau}f (x,y,\tau)|^{2} d\tau \right)^{\frac{1}{2}} \nonumber \\
&&  \leq C\left(\int_{0}^{+\infty}| f(x,y,z)|^{2} dz+\int_{0}^{+\infty} |\partial_{z}f (x,y,z)|^{2} dz \right).
\label{4.29}
\end{eqnarray}
Secondly, we apply the two-dimensional Sobolev inequality for $x, y$ variables to get
\begin{eqnarray}
  |f(x,y,z)|^{2} \leq C(\|f( z)\|_{L^{2}(\mathbb{T}^{2})}^{2}+ \|\partial_{xy}^{1}f( z)\|_{L^{2}(\mathbb{T}^{2})}^{2}
+\|\partial_{xy}^{2}f( z)\|_{L^{2}(\mathbb{T}^{2})}^{2} )
\label{4.30}
\end{eqnarray}
and
\begin{eqnarray}
 |\partial_{z}f (x,y,z)|^{2} \leq C(\partial_{z}\|f( z)\|_{L^{2}(\mathbb{T}^{2})}^{2}+ \|\partial_{z}\partial_{xy}^{1}f( z)\|_{L^{2}(\mathbb{T}^{2})}^{2}
+\|\partial_{z}\partial_{xy}^{2}f( z)\|_{L^{2}(\mathbb{T}^{2})}^{2} ).
\label{4.31}
\end{eqnarray}
At last, substituting (\ref{4.30}) and (\ref{4.31}) into (\ref{4.29}),  we complete the proof of estimate  (\ref{4.12}).

$\hfill$$\Box$

\begin{Lemma} \label{y4.4}
Let the vector field $(u, Ku, w)$ defined on $\mathbb{T}^2 \times \mathbb{R}^+$ satisfy the condition $\partial_x u+ \partial_y (Ku)+ \partial_z w=0$.  If $s\geq 5$ is an integer, $\gamma\geq 1$, $\sigma>\gamma+\frac{1}{2}$ and $\delta\in (0,1)$, then for any $\varphi \in H^{s,\gamma}_{\sigma,\delta}(\mathbb{T}^{2}\times\mathbb{R}^{+})$,
we have the following inequalities: :\\
$(i)$ For $k=0, 1, 2,  \cdots, s $,
 \begin{eqnarray}
\|(1+z)^{\gamma-1} \partial_{xy}^k(u-U)\|_{L^2}
\leq C_{s, \gamma,\sigma,\delta  }\left(\|\varphi\|_{H^{s,\gamma}_{g}}+\|\partial_{xy}^s U\|_{L^2(\mathbb{T}^2)}\right) ;
\label{4.13}
\end{eqnarray}
(ii) For $k=0, 1, 2,  \cdots, s-1 $, if $K \in C^{s+1}(\mathbb{T}^2)$, then
\begin{eqnarray}
\left \|\frac{\partial_{xy}^{k} w +z\tilde{U}}{1+z}  \right\|_{L^{2}}
 \leq C_{s, \gamma,\sigma,\delta  }\left(\|\varphi\|_{H^{s,\gamma}_{g}}+\|\partial_{xy}^s U\|_{L^2(\mathbb{T}^2)}\right) ,
\label{4.14}
\end{eqnarray}
where $\tilde{U}= \partial_x\partial^{k}_{xy} U + \partial_y\partial^{k}_{xy} (KU)$\\
(iii) For all $| \alpha  | \leq s$,
\begin{align}
&\|(1+z)^{\gamma+\alpha_3} D^\alpha \varphi\|_{L^2} \leq \left\{
\begin{aligned}
&C_{s, \gamma,\sigma,\delta  }\left(\|\varphi\|_{H^{s,\gamma}_{g}}+\|\partial_{xy}^s U\|_{L^2(\mathbb{T}^2)}\right) \ \ \mathrm{if}~\alpha =(s_1, s_2, 0) ,\\
&\|\varphi\|_{H^{s,\gamma}_{g}} \ \ \mathrm{if}~\alpha \neq(s_1, s_2, 0).\\
\end{aligned} \label{4.15}
\right.
\end{align}
(iv) For $k=0, 1, 2,  \cdots, s $,
\begin{align}
&\|(1+z)^{\gamma} g_k\|_{{L^2}(\mathbb{T}^2)} \leq \left\{
\begin{aligned}
&C_{s, \gamma,\sigma,\delta  }\left(\|\varphi\|_{H^{s,\gamma}_{g}}+\|\partial_{xy}^s U\|_{L^2(\mathbb{T}^2)}\right) \ \ \mathrm{if}~k \neq s ,\\
&\|\varphi\|_{H^{s,\gamma}_{g}} \ \ \mathrm{if}~k=s.\\
\end{aligned} \label{4.16}
\right.
\end{align}
(v) For $k=0, 1, 2,  \cdots, s-2$,
\begin{eqnarray}
 \|  \partial_{xy}^{k} u  \|_{L^{\infty}} \leq C_{s, \gamma,\sigma,\delta  }\left(\|\varphi\|_{H^{s,\gamma}_{g}}+\|\partial_{xy}^s U\|_{L^2(\mathbb{T}^2)}\right).
\label{4.17}
\end{eqnarray}
(vi) For $k=0, 1, 2,  \cdots, s-3$, if $K \in C^{s+1}(\mathbb{T}^2)$, then
\begin{eqnarray}
 \left\|  \frac{\partial_{xy}^{k} w}{1+z} \right \|_{L^{\infty}} \leq C_{s, \gamma,\sigma,\delta  }\left(\|\varphi\|_{H^{s,\gamma}_{g}}+\|\partial_{xy}^s U\|_{L^2(\mathbb{T}^2)}+1\right).
\label{4.18}
\end{eqnarray}
(vii) For all $| \alpha  | \leq s-3$,
\begin{eqnarray}
\|(1+z)^{\gamma+\alpha_3} D^\alpha \varphi\|_{L^\infty} \leq C_{s, \gamma,\sigma,\delta  }\|\varphi\|_{H^{s,\gamma}_{g}}
\label{4.19}
\end{eqnarray}
where $C_{s, \gamma,\sigma,\delta  }>0$ is a constant and only depends on $s, \gamma,\sigma,\delta  $.
\end{Lemma}
$\mathbf{Proof.}$
First, by the definition of $\| \cdot \|_{H^{s, \gamma-1}}$ and Lemma \ref{y4.2},  we can directly derive the inequality (\ref{4.13}). Using Lemma \ref{y4.1}, (\ref{4.13}) and the fact that $\partial_x u+ \partial_y (Ku)+ \partial_z w=0$, we have
\begin{align}
 \left \|\frac{\partial_{xy}^{k} w +z \tilde{U}}{1+z}  \right\|_{L^{2}}
 &\leq 2 \| \partial_{xy}^{k }(-\partial_x u- \partial_y (Ku))+\tilde{U}\|_{L^{2}  }  \nonumber \\
&\leq 2 \| \partial_{xy}^{k +1}(u-U)\|_{L^{2}  } +C \sum \limits^{k} \limits_{i=0} \|\partial_y\{\partial_{xy}^{i} K \partial_{xy}^{k-i}(u-U)\}\|_{L^2} \nonumber \\
&  \leq C_{s, \gamma,\sigma,\delta  }\left(\|\varphi\|_{H^{s,\gamma}_{g}}+\|\partial_{xy}^s U\|_{L^2(\mathbb{T}^2)}\right).
\end{align}
The proofs of  (\ref{4.15}) and (\ref{4.16}) are obvious, so we omit it here for brevity. Second, for $k=0, 1, 2,  \cdots, s-2$, applying Lemma \ref{y4.3} and (\ref{4.13}), and (\ref{4.15}), we have
\begin{align}
 \|  \partial_{xy}^{k} (u-U)  \|_{L^{\infty}} &\leq C( \|  \partial_{xy}^{k} (u-U)  \|_{L^{2}}+\|  \partial_{xy}^{k+1} (u-U) \|_{L^{2}}+\|  \partial_{xy}^{k+2 } (u-U) \|_{L^{2}}+\|  \partial_{xy}^{k} \varphi \|_{L^{2}}\nonumber\\
&\quad+\|   \partial_{xy}^{k+1}\varphi\|_{L^{2}}+\|  \partial_{xy}^{k+2 } \varphi \|_{L^{2}}) \nonumber  \\
&\leq C_{s, \gamma,\sigma,\delta  }\left(\|\varphi\|_{H^{s,\gamma}_{g}}+\|\partial_{xy}^s U\|_{L^2(\mathbb{T}^2)}\right).
\end{align}
Hence, by the triangle inequality and $\|\partial_{xy}^k U\|_{L^2(\mathbb{T}^2)} \leq \|\partial_{xy}^s U\|_{L^2(\mathbb{T}^2)}$, we can justify (\ref{4.17}).\\
Next, for $k=0, 1, 2,  \cdots, s-3$, using   (\ref{4.13})-(\ref{4.15}) and Lemma \ref{y4.3}, we conclude that
\begin{align}
 \left\|  \frac{\partial_{xy}^{k} w }{1+z}  \right\|_{L^{\infty}}& \leq \left\|  \frac{z\tilde{U} }{1+z}  \right\|_{L^{\infty}}+\left\|  \frac{\partial_{xy}^{k} w+z\tilde{U} }{1+z}  \right\|_{L^{\infty}} \nonumber \\
 &  \leq C \Big(\|   \tilde{U}   \|_{L^{2}(\mathbb{T}^2)}+\|  \partial_{xy} \tilde{U}   \|_{L^{2}(\mathbb{T}^2)}+\|  \partial_{xy}^{2} \tilde{U}   \|_{L^{2}(\mathbb{T}^2)}+\|(1+z)^{-1}(\partial_{xy}^{k} w+z\tilde{U})\|_{L^2} \nonumber \\
&\quad+\|(1+z)^{-1}(\partial_{xy}^{k} \partial_{z}w+\tilde{U})\|_{L^2}+\|(1+z)^{-2}(\partial_{xy}^{k} w+z\tilde{U})\|_{L^2} \nonumber  \\
&\quad+\|(1+z)^{-1}(\partial_{xy}^{k+1} w+z\partial_{xy}\tilde{U})\|_{L^2} +\|(1+z)^{-1}(\partial_{xy}^{k+2} w+z\partial_{xy}^2\tilde{U})\|_{L^2} \nonumber  \\
&\quad+\|(1+z)^{-1}(\partial_{xy}^{k+1} \partial_{z}w+\partial_{xy}\tilde{U})\|_{L^2}+\|(1+z)^{-2}(\partial_{xy}^{k+1} w+z\partial_{xy}\tilde{U})\|_{L^2} \nonumber  \\
&\quad+\|(1+z)^{-1}(\partial_{xy}^{k+2} \partial_{z}w+\partial_{xy}^2\tilde{U})\|_{L^2}+\|(1+z)^{-2}(\partial_{xy}^{k+2} w+z\partial_{xy}^2\tilde{U})\|_{L^2} \Big )
 \nonumber \\
& \leq C_{s, \gamma,\sigma,\delta  } \left(\|\varphi\|_{H^{s,\gamma}_{g}}+\|\partial_{xy}^s U\|_{L^2(\mathbb{T}^2)}+1\right).
\label{4.20}
\end{align}
Last, inequality (\ref{4.19}) follows directly from  inequality (\ref{4.15}) in Lemma \ref{y4.3}.

$\hfill$$\Box$

\section{Uniform estimates }
In this section, we will estimate the norm $\|\varphi\|_{H^{s,\gamma}_{g}}^{2}$ to   problem $(\ref{1.5})$ by the energy method, whose idea  comes from \cite{3}.
\subsection{Estimates on $D^{\alpha}\varphi$ }
We will estimate the norm with weight $(1+z)^{ \gamma+ \alpha_{3}}$ on $D^{\alpha}\varphi$, $\alpha=(\alpha_{1},\alpha_{2}, \alpha_{3})$, $|\alpha|\leq s$, and $\alpha_{1}+\alpha_{2}$ should be restricted as $\alpha_{1}+\alpha_{2}\leq s-1$. Otherwise, it is impossible to directly estimate the norm with $\partial _{xy}^{s}\varphi$.
\begin{Lemma}(Reduction of boundary data)\label{y2.1}
If $\varphi$ solves $(\ref{1.4})$ and $(\ref{1.5})$, then on the boundary at $z=0$, we have,
\begin{equation}\left\{
\begin{array}{ll}
\partial _{z}\varphi|_{z=0}=\partial _{x}P,\\
\partial _{z}^{3}\varphi|_{z=0}=(\partial _{t}-\epsilon^2\partial_x^2-\epsilon^2\partial_y^2)\partial _{x}P+(\varphi \partial_x \varphi+K\varphi \partial_y \varphi)\big|_{z=0}-\partial_y K |\varphi|^2 \big|_{z=0}. \label{3.21}
\end{array}
         \right.\end{equation}
For any $2 \leq k \leq \frac{s}{2}$, there are some constants  $C_{K,k,l,\rho^1,\rho^2,...,\rho^j}$,   not depending on $\epsilon$ (in this section, the superscript  $\epsilon$ is removed for simplification) or $(u,Ku, w)$, such that
\begin{eqnarray}
\begin{aligned}
\partial_z^{2k+1} \varphi|_{z=0}&=(\partial_t-\epsilon^2\partial^2_x-\epsilon^2\partial^2_y)^k\partial_x P \\
&\quad+\sum\limits_{l=0}^{2k-1}\partial_{xy}^{l} K\sum\limits_{l=0}^{k-1}\epsilon^{2l}\sum\limits_{j=2}^{\max\{2,k-l\}}\sum\limits_{\rho \in A^j_{k,l}}C_{K,k,l,\rho^1,\rho^2,...,\rho^j}\prod^{j}_{i=1}D^{\rho^i}\varphi \big{|}_{z=0},
\label{3.22}
\end{aligned}
\end{eqnarray}
where $A^j_{k,l}:=\{\rho:=(\rho^1,\rho^2,...,\rho^j)\in \mathbb{N}^{2j};3\sum\limits_{i=1}^{j}\sum\limits_{i=1}^{j}\rho^i_2=2k+4l+1,\sum\limits_{i=1}^{j}\rho^i_1 \leq k+2l-1,\sum\limits_{i=1}^{j}\rho_2^j\leq 2k-2l-2, ~and~ |\rho^i| \leq 2k-l-1 ~for~ all~ i=1,2,...,j\}$.
\end{Lemma}
$\mathbf{Proof.}$
According to the equation $(\ref{1.5})$ and the boundary condition $\partial _{z}\varphi|_{z=0}=\partial _{x}P$, we get
\begin{align}
\partial _{z}^{3}\varphi\big{|}_{z=0}&= \partial_z\left\{\partial _{t}\varphi+(u\partial _{x} +Ku\partial _{y}+w\partial_{z})\varphi -\varepsilon\partial _{x}^{2}\varphi-\varepsilon\partial _{y}^{2}\varphi- \frac{1}{2}\partial_y K \partial_z |u|^2 \right\}\big|_{z=0}\nonumber\\
&=(\partial _{t}-\epsilon^2\partial_x^2-\epsilon^2\partial_y^2)\partial _{x}P+(\varphi \partial_x \varphi+K\varphi \partial_y \varphi)\big|_{z=0}-\partial_y K |\varphi|^2 \big|_{z=0},
\label{3.23}
\end{align}
and
\begin{eqnarray}
\begin{aligned}
\varphi _{z}^{2n+1}\varphi \big{|}_{z=0}&= \left\{\left(\partial _{t}-\epsilon^2\partial_x^2\varphi-\epsilon^2\partial_y^2\varphi\right)\partial _{z}^{2n-1}\varphi+\partial _{z}^{2n-1}(u\partial _{x} \varphi +Ku\partial _{y}\varphi+w\partial_{z}\varphi-\frac{1}{2}\partial_y K \partial_z |u|^2)\right\}\big|_{z=0}.
\label{3.24}
\end{aligned}
\end{eqnarray}
Hence, for $n=2$, we arrive at
\begin{align}
\partial _{z}^{5}w \big|_{z=0}
&= \left\{(\partial _{t}-\epsilon^2\partial_x^2-\epsilon^2\partial_y^2)\partial_z^{3}\varphi+\partial _{z}^{3}(u\partial _{x} +Ku\partial _{y}+w\partial_{z})\varphi-\partial _{z}^{3}(\frac{1}{2}\partial_y k \partial_z |u|^2) \right\}\Big|_{z=0}\nonumber\\
&=(\partial _{t}-\epsilon^2\partial_x^2-\epsilon^2\partial_y^2)^2\partial_x P+(\partial _{t}-\epsilon^2\partial_x^2-\epsilon^2\partial_y^2)\left\{(\varphi \partial_x \varphi+K\varphi \partial_y \varphi)-\partial_y K |\varphi|^2\right\}\big|_{z=0}\nonumber\\
&\quad+(3\varphi\partial_x\partial_z^2 \varphi+2\partial_z \varphi \partial_x \partial_z \varphi-2\partial_x \varphi \partial^2_z \varphi)|_{z=0}+K(3\varphi\partial_y\partial_z^2 \varphi+2\partial_z \varphi \partial_y \partial_z \varphi-2\partial_y \varphi \partial^2_z \varphi)|_{z=0}\nonumber\\
&\quad-\partial_y K(4\varphi\partial_z^2 \varphi+3\partial_z\varphi\partial_z \varphi)\big|_{z=0}.
\label{3.25}
\end{align}
Since the right-hand side of (\ref{3.25}) except for the second term in the formula is in the desired form, we only need to deal with the terms $(\partial _{t}-\epsilon^2\partial_x^2-\epsilon^2\partial_y^2)\left\{(\varphi \partial_x \varphi+K\varphi \partial_y \varphi)-\partial_y K |\varphi|^2\right\}\big|_{z=0}$. Using the evolution equations for $\varphi$, $\partial_x\varphi$ and $\partial_y\varphi$ as well as $(u, w) \big|_{z=0}=0$, we may check that
\begin{align}
&(\partial _{t}-\epsilon^2\partial_x^2-\epsilon^2\partial_y^2)\left\{(\varphi \partial_x \varphi+K\varphi \partial_y \varphi)-\partial_y K |\varphi|^2\right\}\big|_{z=0}\nonumber \\
&=\left\{\partial_t \varphi (\partial_x \varphi +K\partial_y \varphi -2\partial_y  K \varphi)+\varphi(\partial _{t}\partial_x \varphi+K\partial _{t}\partial_y \varphi)\right\}\big|_{z=0}\nonumber \\
&\quad-(\epsilon^2\partial_x^2+\epsilon^2\partial_y^2)\left\{(\varphi \partial_x \varphi+K\varphi \partial_y \varphi)-\partial_y K |\varphi|^2\right\}\big|_{z=0} \nonumber \\
&= \left( \varepsilon^2\partial _{x}^{2}\varphi+\varepsilon^2\partial _{y}^{2}\varphi+ \partial _{z}^{2}\varphi \right)(\partial_x \varphi +K\partial_y \varphi -2\partial_y  K \varphi)\big|_{z=0} \nonumber \\
&\quad+\varphi \left\{\left( \varepsilon^2\partial _{x}^{3}\varphi+\varepsilon^2\partial _x\partial _{y}^{2}\varphi+ \partial _x\partial _{z}^{2}\varphi \right)+K\left( \varepsilon^2\partial _{y}\partial _{x}^{2}\varphi+\varepsilon^2\partial _{y}^{3}\varphi+ \partial _y\partial _{z}^{2}\varphi \right)\right\}\big|_{z=0} \nonumber \\
& \quad-(\epsilon^2\partial_x^2+\epsilon^2\partial_y^2)\left\{(\varphi \partial_x \varphi+K\varphi \partial_y \varphi)-\partial_y K |\varphi|^2\right\}\big|_{z=0},
\label{3.26}
\end{align}
where  we have used the following facts,
\begin{align*}
\partial _{t} \varphi|_{z=0}&= \left\{-(u\partial _{x} +Ku\partial _{y}+w\partial_{z})\varphi +\varepsilon\partial _{x}^{2}\varphi+\varepsilon\partial _{z}^{2}\varphi+ \partial _{y}^{2}\varphi+\frac{1}{2}\partial_y K \partial_z |u|^2 \right\}\big{|}_{z=0} \nonumber\\
&=\left( \varepsilon^2\partial _{x}^{2}\varphi+\varepsilon^2\partial _{y}^{2}\varphi+ \partial _{z}^{2}\varphi \right)\big{|}_{z=0}, \nonumber\\
\partial _{t}\partial_x \varphi|_{z=0}&= \partial_x\left\{-(u\partial _{x} +Ku\partial _{y}+w\partial_{z})\varphi +\varepsilon\partial _{x}^{2}\varphi+\varepsilon\partial _{z}^{2}\varphi+ \partial _{y}^{2}\varphi+\frac{1}{2}\partial_y K \partial_z |u|^2 \right\}\big{|}_{z=0} \nonumber\\
&=\left( \varepsilon^2\partial _{x}^{3}\varphi+\varepsilon^2\partial _x\partial _{y}^{2}\varphi+ \partial _x\partial _{z}^{2}\varphi \right)\big{|}_{z=0}, \nonumber\\
\partial _{t}\partial_y \varphi|_{z=0}&=\left( \varepsilon^2\partial _{y}\partial _{x}^{2}\varphi+\varepsilon^2\partial _{y}^{3}\varphi+ \partial _y\partial _{z}^{2}\varphi \right)\big{|}_{z=0}.
\end{align*}
We have completed the verification of formula  (\ref{3.22}) for $k=2$.

Now, using the same algorithm, we are going to prove formula (\ref{3.22}) by induction on $k$. For notational convenience, we denote
\begin{eqnarray}
\mathcal{A}_k:=\bigg\{\sum\limits_{l=0}^{2k-1}\partial_{xy}^{l} K\sum\limits_{l=0}^{k-1}\epsilon^{2l}\sum\limits_{j=2}^{\max\{2,k-l\}}\sum\limits_{\rho \in A^j_{k,l}}C_{\wedge_{\alpha},k,l,\rho^1,\rho^2,...,\rho^j}\prod^{j}_{i=1}D^{\rho^i}\varphi\big{|}_{z=0}\bigg\}.
\label{3.27}
\end{eqnarray}
Using this notation, we will prove $\partial_{z}^{2k+1}\varphi\big{|}_{z=0}-(\partial_t-\epsilon^2\partial_x^2-\epsilon^2\partial_y^2)^k\partial_x  P \in \mathcal{A}_k$. Assuming that formula (\ref{3.22}) holds for $k=n$, we will show that it also holds for $k=n+1$ as follows. Then we differentiate the vorticity equation $2n+1$ times  with respect to $z$ to obtain
\begin{align}
\partial _{z}^{2n+3}\varphi\big{|}_{z=0}&= \left\{\left(\partial _{t}-\epsilon^2\partial_x^2\varphi-\epsilon^2\partial_y^2\varphi\right)\partial _{z}^{2n+1}\varphi+\partial _{z}^{2n+1}(u\partial _{x} \varphi +Ku\partial _{y}\varphi+w\partial_{z}\varphi-\partial_y K \varphi u)\right\}\big|_{z=0}\nonumber\\
&=\left(\partial _{t}-\epsilon^2\partial_x^2\varphi-\epsilon^2\partial_y^2\varphi\right)\partial _{z}^{2n+1}\varphi \big|_{z=0} \nonumber\\
&\quad  + \sum^{2n+1}_{m=1}\binom{2n+1}{m} \partial^{m-1}_z \varphi (\partial_x+K\partial_y-\partial_y K) \partial_z^{2n-m+1}\varphi\big|_{z=0} \nonumber\\
&\quad - \sum^{2n+1}_{m=2}\binom{2n+1}{m} (\partial_x+K\partial_y) \partial^{m-2}_z \varphi  \partial_z^{2n-m+2}\varphi\big|_{z=0}.
\label{3.28}
\end{align}
By routine checking, one may show that the last two terms of (\ref{3.28}) belong to
$\mathcal{A}_{k+1}$, so it only remains to deal with the term $(\partial _{t}-\epsilon^2\partial_x^2-\epsilon^2\partial_y^2)\partial _{z}^{2n+1}w$.
Thanks to the induction hypothesis, there exist constants $C_{K,n,l,\rho^1,\rho^2,...,\rho^j}$ such that
\begin{align}
\partial_z^{2n+1}\varphi\big{|}_{z=0}&=(\partial_t-\epsilon^2\partial^2_x-\epsilon^2\partial^2_y)^n\partial_x P\nonumber\\
&\quad+\sum\limits_{l=0}^{2k-1}\partial_{xy}^{l}K\sum\limits_{l=0}^{k-1}\epsilon^{2l}\sum\limits_{j=2}^{\max\{2,n-l\}}\sum\limits_{\rho \in A^j_{n,l}}C_{K,n,l,\rho^1,\rho^2,...,\rho^j}\prod^{j}_{i=1}D^{\rho^i}\varphi\big{|}_{z=0}.
\label{3.29}
\end{align}
Thus we have, up to a relabeling of the indices $\rho^i$,
\begin{align}
&(\partial_t-\epsilon^2\partial^2_x-\epsilon^2\partial^2_y)\partial_z^{2n+1}\varphi\big{|}_{z=0} \nonumber\\
&=(\partial_t-\epsilon^2\partial^2_x-\epsilon^2\partial^2_y)^{n+1}\partial_x P \nonumber\\
&\quad+\sum\limits_{l=0}^{2k-1}\partial_{xy}^{l}K\sum\limits_{l=0}^{n-1}\epsilon^{2l}\sum\limits_{j=2}^{\max\{2,n-l\}}\sum\limits_{\rho \in A^j_{n,l}} \widetilde{C}_{K,n,l,\rho^1,\rho^2,...,\rho^j}(\partial_t-\epsilon^2\partial^2_x-\epsilon^2\partial^2_y)D^{\rho^1}\varphi\prod^{j}_{i=2}D^{\rho^i}\varphi\big{|}_{z=0}\nonumber\\
&\quad-\sum\limits_{l=0}^{2k-1}\partial_{xy}^{l+2}K\sum\limits_{l=0}^{n-1}\epsilon^{2l+2}\sum\limits_{j=2}^{\max\{2,n-l\}}\sum\limits_{\rho \in A^j_{n,l}}\widetilde{\widetilde{C}}_{K,n,l,\rho^1,\rho^2,...,\rho^j}\partial_x D^{\rho^1}\varphi\partial_x D^{\rho^2}\varphi\prod^{j}_{i=3}D^{\rho^i}\varphi\big{|}_{z=0},
\label{3.30}
\end{align}
where $\widetilde{C}_{K,n,l,\rho^1,\rho^2,...,\rho^j}$ and $\widetilde{\widetilde{C}}_{K,n,l,\rho^1,\rho^2,...,\rho^j}$ are some new constants depending on $C_{K,n,l,\rho^1,\rho^2,...,\rho^j}$. It is worth noting that the last term on the right-hand side of (\ref{3.30})
belongs to $\mathcal{A}_{n+1}$, so it remains to check whether the second term on right-hand
side of (\ref{3.30}) also belongs to $\mathcal{A}_{n+1}$.

Indeed, differentiating the vorticity equation $(\ref{1.5})_1$ $\rho_1^1$ times,   $\rho_2^1$ times, $\rho_3^1$ times  with respect to $x,y,z$ respectively,  and then evaluating at $z=0$, we have, by denoting $e_3=(0,0,1)$,
\begin{align}
&(\partial_t-\epsilon^2\partial_x^2-\epsilon^2\partial_y^2)D^{\rho^1}\varphi \big{|}_{z=0} \nonumber \\
&=-\sum_{\substack{\beta \le \rho^1 \\
\beta_3 \geq 1}}\binom{\rho^1}{\beta} D^{\beta-e_3} \varphi \partial_x D^{\rho_1-\beta} \varphi \big{|}_{z=0}
-\sum_{\substack{\beta \le \rho^1 \\
\beta_3 \geq 1}}\binom{\rho^1}{\beta} D^{\beta-e_3} (K\varphi) \partial_y D^{\rho_1-\beta} \varphi \big{|}_{z=0}\nonumber\\
&\quad +\sum_{\substack{\beta \le \rho^1 \\
\beta_3 \geq 2}}\binom{\rho^1}{\beta} D^{\beta-2e_3} (\partial_x \varphi +\partial_y(K\varphi)) \partial_z D^{\rho_1-\beta} w \big{|}_{z=0}+\partial_z^2 D^{\rho^1}\varphi \big{|}_{z=0}.
\label{3.31}
\end{align}
Using (\ref{3.31}), one may justify by a routine counting of indices that the second term
on the right-hand side of (\ref{3.30}) belongs to $\mathcal{A}_{n+1}$.
This thus completes the proof.
$\hfill$$\Box$

\begin{Lemma}\label{y2.2}
Let $s \geq 5$ be an even integer, $\gamma\geq 1$, $\sigma>\gamma+\frac{1}{2}$ and $\delta\in (0,1)$. If  $\varphi \in H^{s,\gamma}_{\sigma,\delta}$  solves $(\ref{1.4})$  and $(\ref{1.5})$, we have the following  estimates:\\
(i)\;when $|\alpha|\leq s-1$,
\begin{eqnarray}
\bigg|\iint_{\mathbb{T}^2 } D^{\alpha}\varphi\partial _{z}D^{\alpha} \varphi dxdy \big{|} _{z=0}\bigg|\leq \frac{1}{12}\|(1+z)^{\gamma+\alpha_3+1}\partial^2_z D^\alpha \varphi\|_{L^2}^2+C\|\varphi\|^2_{H^{s,\gamma}_g},
\end{eqnarray}
(ii)\;when $|\alpha|=s$, $\alpha_3$ is even,
\begin{align}
\bigg|\iint_{\mathbb{T} ^2} D^{\alpha}\varphi\partial _{z}D^{\alpha}\varphi dxdy \big{|} _{z=0}\bigg|\leq &\frac{1}{12}\|(1+z)^{\gamma+\alpha_2}\partial_z D^\alpha \varphi\|^2_{L^2}+ \sum\limits_{l=0}^{\frac{s}{2}}\|\partial^{l}_t\partial_xP\|_{H^{s-2l}(\mathbb{T}^2)}^{2}\nonumber\\
&+C_{s,\gamma,\sigma,\delta}\|\partial_{xy}^s K\|_{L^{\infty}(\mathbb{T}^2)}(1+\|\varphi\|^2_{H^{s,\gamma}_g})^{s-2}\|\varphi\|^2_{H^{s,\gamma}_g},
\end{align}
(iii)\;when $|\alpha|=s$, $\alpha_3$ is odd,
\begin{align}
\bigg|\iint_{\mathbb{T} ^2} D^{\alpha}\varphi\partial _{z}D^{\alpha}\varphi dxdy \big{|} _{z=0}\bigg|\leq& \frac{1}{12}\|(1+z)^{\gamma+\alpha_3+1}\partial_x^{\alpha_1-1}\partial_y^{\alpha_2}\partial_z^{\alpha_3+2} \varphi\|^2_{L^2}+\sum\limits_{l=0}^{\frac{s}{2}}\|\partial^{l}_t\partial_xP\|_{H^{s-2l}(\mathbb{T})^2}^{2}\nonumber\\
&+C_{s,\gamma,\sigma,\delta}\|\partial_{xy}^s K\|_{L^{\infty}(\mathbb{T}^2)}(1+\|\varphi\|^2_{H^{s,\gamma}_g})^{s-2}\|\varphi \|^2_{H^{s,\gamma}_g}.
\end{align}
\end{Lemma}
$\mathbf{Proof.}$
\emph{Case 1.} When $|\alpha|\leq s-1$,  using the following trace estimate
\begin{eqnarray}
\iint_{\mathbb{T}^2 } |f|dxdy \big{|}_{z=0} \leq C\left(\int^1_0\int_{\mathbb{T}^2}|f|dxdydz+\int^1_0\int_{\mathbb{T}^2}|\partial_z f|dxdydz\right),
\end{eqnarray}
we have
\begin{eqnarray}
\bigg|\int_{\mathbb{T}^2 } D^{\alpha}\varphi \partial _{z}D^{\alpha}\varphi dxdy \big{|} _{z=0}\bigg|\leq \frac{1}{12}\|(1+z)^{\gamma+\alpha_3+1}\partial^2_z D^\alpha \varphi\|_{L^2}^2+C\|\varphi\|^2_{H^{s,\gamma}_g}.
\end{eqnarray}
\emph{Case 2.} When $|\alpha|=s$, $\alpha_3=2k$ for some $k\in \mathbb{N}$, we can apply   Lemma \ref{y2.1} to $\partial_z D^\alpha \varphi | _{z=0}$ to obtain
\begin{align}
&\int_{\mathbb{T} } D^{\alpha}\varphi\partial _{z}D^{\alpha}\varphi dxdy \big{|} _{z=0} \nonumber \\
&=\int_{\mathbb{T}^2 } D^{\alpha}\varphi(\partial_t-\epsilon^2\partial_x^2-\epsilon^2\partial_y^2)^k\partial_x^{\alpha_1}\partial_y^{\alpha_2}\partial_x P dx dy\big{|} _{z=0}\nonumber\\
&\quad+\sum\limits_{l=0}^{2k-1}\partial_{xy}^{l} K\sum\limits_{l=0}^{k-1}\epsilon^{2l}\sum\limits_{j=2}^{\max\{2,k-l\}}\sum\limits_{\rho \in A^j_{k,l}}C_{K,k,l,\rho^1,\rho^2,...,\rho^j}\int_{\mathbb{T}^2 }D^\alpha \varphi \partial_x^{\alpha_1}\partial_y^{\alpha_2}(\prod^{j}_{i=1}D^{\rho^i}\varphi)dx\big{|}_{z=0}.
\end{align}
Then we again apply the simple trace estimate to control the boundary integral as follows
\begin{align}
\bigg|\int_{\mathbb{T} ^2} D^{\alpha}\varphi \partial _{z}D^{\alpha}\varphi dxdy \big{|} _{z=0}\bigg|\leq& \frac{1}{12}\|(1+z)^{\gamma+\alpha_3}\partial_z D^\alpha \varphi\|^2_{L^2}+ \sum\limits_{l=0}^{\frac{s}{2}}\|\partial^{l}_t\partial_xP\|_{H^{s-2l}(\mathbb{T}^2)}^{2}\nonumber\\
&+C_{s,\gamma,\sigma,\delta}\|\partial_{xy}^s K\|_{L^{\infty}(\mathbb{T}^2)}(1+\|\varphi\|_{H_g^{s,\gamma}})^{s-2}\|\varphi\|^2_{H^{s,\gamma}_g}.
\end{align}
\emph{Case 3.} When $|\alpha|=s$, $\alpha_3=2k+1$ for some $k\in \mathbb{N}$, using integration by parts in the $x,y$ variables, we have
\begin{eqnarray}
\int_{\mathbb{T}^2 } D^{\alpha}\varphi \partial _{z}D^{\alpha}\varphi dxdy \big{|} _{z=0}=-\int_{\mathbb{T}^2}\partial_x D^\alpha \varphi \partial_x^{\alpha_1-1}\partial_y^{\alpha_2}\partial_z^{\alpha_3+1}\varphi dxdy \big{|}_{z=0}.
\end{eqnarray}
The term $\partial_x D^\alpha \varphi \big{|}_{z=0}=\partial_x^{\alpha_1-1}\partial_y^{\alpha_2}\partial_z^{\alpha_3+1}\varphi \big{|} _{z=0}$ has an odd number of derivatives in $z$, then we get
\begin{align}
\bigg|\int_{\mathbb{T} ^2} D^{\alpha}\varphi\partial _{y}D^{\alpha}\varphi dxdy \big{|} _{z=0}\bigg|\leq& \frac{1}{12}\|(1+z)^{\gamma+\alpha_3+1}\partial_x^{\alpha_1-1}\partial_y^{\alpha_2}\partial_z^{\alpha_3+2}\varphi \|^2_{L^2}+ \sum\limits_{l=0}^{\frac{s}{2}}\|\partial^{l}_t\partial_xP\|_{H^{s-2l}(\mathbb{T}^2)}^{2}\nonumber\\
&+C_{s,\gamma,\sigma,\delta}\|\partial_{xy}^s K\|_{L^{\infty}(\mathbb{T}^2)}(1+\|\varphi\|_{H_g^{s,\gamma}})^{s-2}\|\varphi\|^2_{H^{s,\gamma}_g}.
\end{align}
$\hfill$$\Box$

\begin{Proposition}\label{p2.1}
Let $s \geq 5$ be an even integer, $\gamma\geq 1$, $\sigma>\gamma+\frac{1}{2}$ and $\delta\in (0,1)$. If  $\varphi \in H^{s,\gamma}_{\sigma,\delta}$  solves $(\ref{1.4})$  and $(\ref{1.5})$, then we have the following uniform (in $\epsilon$) estimate:
\begin{eqnarray}
\begin{aligned}
\frac{d}{dt}\|(1+z)^{ \gamma+ \alpha_{3}}D^{\alpha}\varphi\|_{L^{2}}^{2} &\leq C_{s, \gamma,\sigma,\delta  }\left(1+\|\partial_{xy}^{s} K\|_{L^\infty(\mathbb{T}^2)}\right)\left(\|\varphi\|_{H^{s,\gamma}_{g}}+\|\partial_{xy}^s U\|_{L^2(\mathbb{T}^2)}+1\right)\|\varphi \|_{H^{s,\gamma}_{g}}^{2} \\
& \quad +C_{s,\gamma,\sigma,\delta}\|\partial_{xy}^s K\|_{L^{\infty}(\mathbb{T}^2)}(1+\|w\|_{H_g^{s,\gamma}})^{s-2}\|\varphi\|^2_{H^{s,\gamma}_g}\\
&\quad +\sum\limits_{l=0}^{\frac{s}{2}}\|\partial^{l}_t\partial_xP\|_{H^{s-2l}(\mathbb{T}^2)}^{2}. \label{2.13}
\end{aligned}
\end{eqnarray}
\end{Proposition}
$\mathbf{Proof.}$
Applying the operator $D^{\alpha}$ on $(\ref{1.5})_{1}$ with $\alpha=(\alpha_{1},\alpha_{2},\alpha_{3}),~|\alpha|\leq s, ~\alpha_{1}+\alpha_{2}\leq s-1$,
\begin{align}
&(\partial _{t} +u\partial _{x} +Ku\partial _{y}+w\partial _{z}-\varepsilon^2\partial _{x}^{2}\varphi-\varepsilon^2\partial _{y}^{2}\varphi-\partial _{z}^{2})D^{\alpha}\varphi \nonumber\\
&=- \sum\limits_{0<\beta\leq \alpha}\binom{ \alpha}{\beta} \left\{D^{\beta} u\partial _{x}D^{\alpha-\beta}\varphi  + D^{\beta}(Ku)\partial _{y}D^{\alpha-\beta}\varphi + D^{\beta} w\partial _{z}D^{\alpha-\beta}\varphi \right\}\nonumber\\
&\quad+ \sum\limits_{0 \leq \beta\leq \alpha}\binom{ \alpha}{\beta} D^{\beta} (u\partial_y K )D^{\alpha-\beta}\varphi .
 \label{2.1}
 \end{align}
Multiplying $(\ref{2.1})$ by $(1+z)^{2\gamma+2\alpha_{3}}D^{\alpha}\varphi$ and integrating it over $\mathbb{T}^{2}\times \mathbb{R}^{+}$ yield
\begin{eqnarray}
\begin{aligned}
&\frac{1}{2}\frac{d}{dt}\|(1+z)^{ \gamma+ \alpha_{3}}D^{\alpha}\varphi\|_{L^{2}}^{2}
+\|(1+z)^{ \gamma+ \alpha_{3}}\partial _{z}D^{\alpha}\varphi\|_{L^{2}}^{2}+\varepsilon^2\|(1+z)^{ \gamma+ \alpha_{3}}\partial _{z}D^{\alpha}\varphi\|_{L^{2}}^{2}+\varepsilon^2\|(1+z)^{ \gamma+ \alpha_{3}}\partial _{z}D^{\alpha}\varphi\|_{L^{2}}^{2}   \\
&=-\iint_{\mathbb{T} ^{2}} D^{\alpha}\varphi\partial _{z}D^{\alpha}\varphi dxdy | _{z=0}-(2\gamma+2\alpha_{3}) \iiint(1+z)^{2\gamma+2\alpha_{3}-1} D^{\alpha}\varphi\partial _{z}   D^{\alpha}\varphi \\
&\quad +(\gamma+\alpha_{3}) \iiint(1+z)^{2\gamma+2\alpha_{3}-1} w| D^{\alpha}\varphi|^{2}  \\
&\quad  - \sum\limits_{0<\beta\leq \alpha}\binom{ \alpha}{\beta} \iiint(1+z)^{2\gamma+2\alpha_{3}}D^{\alpha}\varphi  \left\{D^{\beta} u\partial _{x}D^{\alpha-\beta}\varphi + D^{\beta}(Ku)\partial _{y}D^{\alpha-\beta}\varphi+ D^{\beta} w\partial _{z}D^{\alpha-\beta}\varphi \right\} \\
&\quad  + \sum\limits_{0<\beta\leq \alpha}\binom{ \alpha}{\beta} \iiint(1+z)^{2\gamma+2\alpha_{3}}D^{\alpha}\varphi  \left\{D^{\beta} (u \partial_y K )D^{\alpha-\beta}\varphi  \right\}.
\label{2.2}
\end{aligned}
\end{eqnarray}
Indeed, $(\ref{2.2})$ can be obtained by integration by parts and by using the boundary condition, i.e.,
\begin{eqnarray}
\begin{aligned}
& \iiint(1+z)^{2\gamma+2\alpha_{3}}D^{\alpha} \varphi(u\partial _{x} +ku\partial _{y} +w\partial _{z} )D^{\alpha}\varphi \\
&=-(\gamma +\alpha_{3})\iiint(1+z)^{2\gamma+2\alpha_{3}-1} w |D^{\alpha}\varphi |^{2},
 \label{2.3}
\end{aligned}
\end{eqnarray}
and
\begin{eqnarray}
\begin{aligned}
 -\iiint(1+z)^{2\gamma+2\alpha_{3}}D^{\alpha}\varphi \partial _{z}^{2}D^{\alpha}\varphi & = \|(1+z)^{\gamma+ \alpha_{3}}\partial _{z}D^{\alpha}\varphi\|_{L^{2}}+\iint_{\mathbb{T}^{2} } D^{\alpha}\varphi\partial _{z}D^{\alpha}\varphi dxdy | _{z=0}\\
& \quad +  (2\gamma+2\alpha_{3})  \iiint(1+z)^{2\gamma+2\alpha_{3}-1}D^{\alpha}\varphi\partial _{z} D^{\alpha}\varphi.
 \label{2.4}
\end{aligned}
\end{eqnarray}
We need to estimate   equation $(\ref{2.2})$ term by term. Obviously, the first term on the right-hand side of $(\ref{2.2})$ follows from Lemmas \ref{y2.1}-\ref{y2.2}.

Secondly, using the H$\ddot{o}$lder inequality, we have
\begin{eqnarray}
&&\left|(2\gamma+2\alpha_{3}) \iiint(1+z)^{2\gamma+2\alpha_{3}-1} D^{\alpha}\varphi\partial _{z}   D^{\alpha}\varphi \right|\nonumber \\
 &\leq& (2\gamma+2\alpha_{3}) \|\frac{1}{1+z }  \|_{L^{\infty}} \|(1+z)^{\gamma+ \alpha_{3}}D^{\alpha}\varphi\|_{L^{2}}
  \|(1+z)^{\gamma+ \alpha_{3}}\partial_{z} D^{\alpha}\varphi \|_{L^{2}} \nonumber \\
  & \leq & C_{ \gamma,\alpha } \|(1+z)^{\gamma+ \alpha_{3}}D^{\alpha}\varphi \|_{L^{2}}^{2}+\frac{1}{4}\|(1+z)^{\gamma+ \alpha_{3}}\partial_{z} D^{\alpha}\varphi \|_{L^{2}}^{2} \nonumber \\
& \leq & C_{ \gamma,\alpha } \|\varphi \|_{H^{s, \gamma}_g}^{2}+\frac{1}{4}\|(1+z)^{\gamma+ \alpha_{3}}\partial_{z} D^{\alpha}\varphi \|_{L^{2}}^{2} ,  \label{2.5}
\end{eqnarray}
which, together with Lemma \ref{y4.4}, yields
\begin{eqnarray}
\left|   (\gamma+\alpha_{3}) \iiint(1+z)^{2\gamma+2\alpha_{3}-1} w  | D^{\alpha}\varphi|^{2}  \right|
 &\leq&  \|\frac{w}{1+z}\|_{L^{\infty}}  \|(1+z)^{\gamma+ \alpha_{3}}D^{\alpha}\varphi\|_{L^{2}}^{2} \nonumber \\
& \leq &  C_{s, \gamma,\sigma,\delta  }\left(\|\varphi\|_{H^{s,\gamma}_{g}}+\|\partial_{xy}^s U\|_{L^2(\mathbb{T}^2)}+1\right)\|\varphi\|_{H^{s,\gamma}_{g}}^2.  \label{2.6}
\end{eqnarray}
Last,  it follows that the first term of the third line on the right-hand side of equation $(\ref{2.2})$ has the following three cases,
\begin{eqnarray}
J_1=- \sum\limits_{0<\beta\leq \alpha}\binom{ \alpha}{\beta} \iiint(1+z)^{2\gamma+2\alpha_{3}}D^{\alpha}\varphi  D^{\beta} u\partial _{x}D^{\alpha-\beta}\varphi,
\end{eqnarray}
\begin{eqnarray}
J_2=- \sum\limits_{0<\beta\leq \alpha}\binom{ \alpha}{\beta} \iiint(1+z)^{2\gamma+2\alpha_{3}}D^{\alpha}\varphi  D^{\beta}(Ku)\partial _{y}D^{\alpha-\beta}\varphi,
\end{eqnarray}
and
\begin{eqnarray}
J_3=- \sum\limits_{0<\beta\leq \alpha}\binom{ \alpha}{\beta} \iiint(1+z)^{2\gamma+2\alpha_{3}}D^{\alpha}\varphi   D^{\beta} w\partial _{z}D^{\alpha-\beta}\varphi.
\end{eqnarray}
Here $J_1,J_2$ and $J_3$ can be estimated in the following.

When $\beta_3 > 0$, we have
\begin{align}
|J_{1}|&\leq   \left|\sum\limits_{0<\beta\leq \alpha}\binom{ \alpha}{\beta} \iiint(1+z)^{2\gamma+2\alpha_{3}}D^{\alpha}\varphi  D^{\beta-e_3} \varphi D^{\alpha-\beta+e_1}\varphi\right| \nonumber \\
& \leq C \| (1+z)^{\gamma+\beta_{3}-1}  D^{\beta-e_3} \cdot (1+z)^{\gamma+\alpha_{3}-\beta_{3}}  D^{\alpha-\beta+e_1}\varphi\|_{L^2} \|(1+z)^{\gamma+\alpha_{3}}D^{\alpha}\|_{L^2}\nonumber \\
& \leq C\|\varphi \|_{H^{s,\gamma}_{g}}^{3}.
\label{2.7}
\end{align}
When $\beta_3 = 0$,  $\beta_1+ \beta_2 = s-1$, which means $\alpha_3 = 0$ or 1, we get
\begin{align}
|J_{1}|&\leq   \left|\sum\limits_{0<\beta\leq \alpha}\binom{ \alpha}{\beta} \iiint(1+z)^{2\gamma+2\alpha_{3}}D^{\alpha}\varphi  \partial_{xy}^{s-1} u \partial_z^{\alpha_3}\varphi\right| \nonumber \\
& \leq C \left(\|  \partial_{xy}^{s-1} (u-U)\|_{L^2}+ \| \partial_{xy}^{s-1} U \|_{L^2}\right) \|(1+z)^{\gamma+\alpha_{3}} \partial_z^{\alpha_3}\varphi\|_{L^\infty}\|(1+z)^{\gamma+\alpha_{3}}D^{\alpha}\|_{L^2}\nonumber \\
& \leq C\left(\left\|\frac{1}{(1+z)^{\gamma-1}}\right\|_{L^\infty}\| (1+z)^{\gamma-1}  \partial_{xy}^{s-1} (u-U)\|_{L^2}+ \| \partial_{xy}^{s-1} U \|_{L^2}\right)\|\varphi \|_{H^{s,\gamma}_{g}}^{2} \nonumber \\
& \leq C_{s, \gamma,\sigma,\delta  }\left(\|\varphi\|_{H^{s,\gamma}_{g}}+\|\partial_{xy}^s U\|_{L^2(\mathbb{T}^2)}\right)\|\varphi \|_{H^{s,\gamma}_{g}}^{2}.
\end{align}
When $\beta_3 = 0$,  $\beta_1+ \beta_2 \leq s-2$, we derive
\begin{align}
|J_{1}|
& \leq C \|  \partial_{x}^{\beta_1}\partial_{y}^{\beta_2} u\|_{L^\infty} \|(1+z)^{\gamma+\alpha_{3}} \partial_{x}^{\alpha_1-\beta_1+1}\partial_{x}^{\alpha_2-\beta_2}\partial_{z}^{\alpha_3}\varphi\|_{L^2}\|(1+z)^{\gamma+\alpha_{3}}D^{\alpha}\|_{L^2}\nonumber \\
& \leq C_{s, \gamma,\sigma,\delta  }\left(\|\varphi\|_{H^{s,\gamma}_{g}}+\|\partial_{xy}^s U\|_{L^2(\mathbb{T}^2)}\right)\|\varphi \|_{H^{s,\gamma}_{g}}^{2}.
\end{align}
Similarly,
when $\beta_3 > 0$, we obtain
\begin{align}
|J_{2}|\leq    C\|\partial_{xy}^s K\|_{L^{\infty}(\mathbb{T}^2)}\|\varphi \|_{H^{s,\gamma}_{g}}^{3}.
\end{align}
When $\beta_3 = 0$,  $\beta_1+ \beta_2 = s-1$,  we conclude
\begin{align}
|J_{2}| \leq C_{s, \gamma,\sigma,\delta  }\|\partial_{xy}^s K\|_{L^{\infty}(\mathbb{T}^2)}\left(\|\varphi\|_{H^{s,\gamma}_{g}}+\|\partial_{xy}^s U\|_{L^2(\mathbb{T}^2)}\right)\|\varphi \|_{H^{s,\gamma}_{g}}^{2}.
\end{align}
When $\beta_3 = 0$,  $\beta_1+ \beta_2 \leq s-2$, we get
\begin{align}
|J_{2}| \leq C_{s, \gamma,\sigma,\delta  }\|\partial_{xy}^s K\|_{L^{\infty}(\mathbb{T}^2)}\left(\|\varphi\|_{H^{s,\gamma}_{g}}+\|\partial_{xy}^s U\|_{L^2(\mathbb{T}^2)}\right)\|\varphi \|_{H^{s,\gamma}_{g}}^{2}.
\end{align}
When $\beta_3 =0$, $\beta_1+ \beta_2 =s-2$ or  $s-1$, we achieve
\begin{align}
|J_{3}|&\leq\left|\sum\limits_{0<\beta\leq \alpha}\binom{ \alpha}{\beta} \iiint(1+z)^{2\gamma+2\alpha_{3}}D^{\alpha}\varphi   \partial_x^{\beta_1}\partial_y^{\beta_2} w\partial _{z}\partial_x^{\alpha_1-\beta_1}\partial_{y}^{\alpha_2 - \beta_2}\partial_z^{\alpha_3+1}\varphi\right| \nonumber \\
& \leq C\|(1+z)^{\gamma+\alpha_{3}}D^{\alpha}\varphi\|_{L^2}\|(1+z)^{-1}(\partial_{x}^{\beta_1}\partial_{y}^{\beta_2} w + z\tilde{U})\|_{L^2}   \|(1+z)^{\gamma+\alpha_3+1} \partial_x^{\alpha_1-\beta_1}\partial_y^{\alpha_2 - \beta_2} \partial_z^{\alpha_3+1} \varphi\|_{L^\infty} \nonumber \\
& \quad +C\|(1+z)^{\gamma+\alpha_{3}}D^{\alpha}\varphi\|_{L^2}\|(1+z)^{-1} z\tilde{U}\|_{L^\infty}
\|(1+z)^{\gamma+\alpha_3+1} \partial_x^{\alpha_1-\beta_1}\partial_y^{\alpha_2 - \beta_2} \partial_z^{\alpha_3+1} \varphi\|_{L^2} \nonumber \\
&  \leq C_{s, \gamma,\sigma,\delta  }\left(\|\varphi\|_{H^{s,\gamma}_{g}}+\|\partial_{xy}^s U\|_{L^2(\mathbb{T}^2)}\right)\|\varphi \|_{H^{s,\gamma}_{g}}^{2}.
\end{align}
When $\beta_3 =0$, $\beta_1+ \beta_2 \leq s-3$, we deduce
\begin{align}
|J_{3}|&\leq\left|\sum\limits_{0<\beta\leq \alpha}\binom{ \alpha}{\beta} \iiint(1+z)^{2\gamma+2\alpha_{3}}D^{\alpha}\varphi   \partial_x^{\beta_1}\partial_y^{\beta_2} w\partial _{z}\partial_x^{\alpha_1-\beta_1}\partial_{y}^{\alpha_2 - \beta_2}\partial_z^{\alpha_3+1}\varphi\right| \nonumber \\
& \leq C\|(1+z)^{\gamma+\alpha_{3}}D^{\alpha}\varphi\|_{L^2}\|(1+z)^{-1}\partial_x^{\beta_1}\partial_y^{\beta_2} w\|_{L^\infty}\|(1+z)^{\gamma+\alpha_3+1} \partial_x^{\alpha_1-\beta_1}\partial_y^{\alpha_2 - \beta_2} \partial_z^{\alpha_3+1} \varphi\|_{L^2}\nonumber \\
& \leq C_{s, \gamma,\sigma,\delta  }\left(\|\varphi\|_{H^{s,\gamma}_{g}}+\|\partial_{xy}^s U\|_{L^2(\mathbb{T}^2)}+1\right)\|\varphi \|_{H^{s,\gamma}_{g}}^{2}.
\end{align}
When $\beta_3 =1$, we obtain
\begin{align}
|J_{3}|&\leq\left|\sum\limits_{0<\beta\leq \alpha}\binom{ \alpha}{\beta} \iiint(1+z)^{2\gamma+2\alpha_{3}}D^{\alpha}\varphi   \partial_x^{\beta_1}\partial_y^{\beta_2} [\partial_x u +\partial_y(Ku)]\partial _{z}\partial_x^{\alpha_1-\beta_1}\partial_{y}^{\alpha_2 - \beta_2}\partial_z^{\alpha_3+1}\varphi\right| \nonumber \\
&  \leq C_{s, \gamma,\sigma,\delta  }(1+\|\partial_{y}\partial_{xy}^{s-1} K\|_{L^\infty(\mathbb{T}^2)})\left(\|\varphi\|_{H^{s,\gamma}_{g}}+\|\partial_{xy}^s U\|_{L^2(\mathbb{T}^2)}+1\right)\|\varphi \|_{H^{s,\gamma}_{g}}^{2}.
\end{align}
When $\beta_3  \geq 2$, we get
\begin{align}
|J_{3}|&\leq\left|\sum\limits_{0<\beta\leq \alpha}\binom{ \alpha}{\beta} \iiint(1+z)^{2\gamma+2\alpha_{3}}D^{\alpha}\varphi   \partial_x^{\beta_1}\partial_y^{\beta_2}\partial_z^{\beta_3-2} [\partial_x \varphi +\partial_yK\varphi +K \partial_y \varphi]\partial _{z}\partial_x^{\alpha_1-\beta_1}\partial_{y}^{\alpha_2 - \beta_2}\partial_z^{\alpha_3+1}\varphi\right| \nonumber \\
&  \leq C_{s, \gamma,\sigma,\delta  }(1+\|\partial_{y}\partial_{xy}^{s-1} K\|_{L^\infty(\mathbb{T}^2)})\|\varphi \|_{H^{s,\gamma}_{g}}^{3}.
\end{align}
The estimate of $J_{4}$ is similar to that of $J_1$, then, when $\beta = 0$, we conclude
\begin{eqnarray}
|J_4| \leq C_{s, \gamma,\sigma,\delta  }\|\partial_{y} K\|_{L^\infty(\mathbb{T}^2)})\left(\|\varphi\|_{H^{s,\gamma}_{g}}+\|\partial_{xy}^s U\|_{L^2(\mathbb{T}^2)}+1\right)\|\varphi \|_{H^{s,\gamma}_{g}}^{2}.
\end{eqnarray}
When $\beta \neq 0$, we achieve
\begin{eqnarray}
|J_4| \leq C_{s, \gamma,\sigma,\delta  }\|\partial_{y}\partial_{xy}^{s-1} K\|_{L^\infty(\mathbb{T}^2)}\left(\|\varphi\|_{H^{s,\gamma}_{g}}+\|\partial_{xy}^s U\|_{L^2(\mathbb{T}^2)}+1\right)\|\varphi \|_{H^{s,\gamma}_{g}}^{2}.
\label{2.12}
\end{eqnarray}
Hence, combining $(\ref{2.2})$, $(\ref{2.5})$-$(\ref{2.6})$ with $(\ref{2.7})$-$(\ref{2.12})$, we conclude
\begin{eqnarray}
\begin{aligned}
\frac{d}{dt}\|(1+z)^{ \gamma+ \alpha_{3}}D^{\alpha}\varphi\|_{L^{2}}^{2} &\leq C_{s, \gamma,\sigma,\delta  }\left(1+\|\partial_{xy}^{s} K\|_{L^\infty(\mathbb{T}^2)}\right)\left(\|\varphi\|_{H^{s,\gamma}_{g}}+\|\partial_{xy}^s U\|_{L^2(\mathbb{T}^2)}+1\right)\|\varphi \|_{H^{s,\gamma}_{g}}^{2} \\
& \quad +C_{s,\gamma,\sigma,\delta}\|\partial_{xy}^s K\|_{L^{\infty}(\mathbb{T}^2)}(1+\|w\|_{H_g^{s,\gamma}})^{s-2}\|\varphi\|^2_{H^{s,\gamma}_g}\\
&\quad +\sum\limits_{l=0}^{\frac{s}{2}}\|\partial^{l}_t\partial_xP\|_{H^{s-2l}(\mathbb{T}^2)}^{2}.
\end{aligned}
\end{eqnarray}
$\hfill$$\Box$

\subsection{Estimates on $g_{s}$ }
In this subsection, we will estimate the norm of $\partial_{x}^{i}\partial_{y}^{s-i}\varphi$ with weight  $(1+z)^{\gamma }$. To overcome the difficult term  $\partial_{z}^{s}w$, we use the cancellation property as usual, i.e., we estimate the norm of $(1+z)^{\gamma}g_{s}$.
\begin{Proposition}\label{p2.2}
Under the same assumption of Proposition \ref{p2.1}, we have the following estimate:
\begin{eqnarray}
\begin{aligned}
\frac{d}{dt}\|(1+z)^{ \gamma}g_{s}\|_{L^{2}}^{2}& \leq  C_{s,\gamma, \sigma, \delta} \left(1+\|\partial_{xy}^s K \|_{L^\infty}\right)\left(1+\|\varphi\|_{H^{s,\gamma}_{g}}+\|\partial_{xy}^s U\|_{L^\infty(\mathbb{T}^2)}\right)\\
&\quad \times\left(\|\varphi\|_{H^{s,\gamma}_{g}}+\|\partial_{xy}^{s+1} U\|_{L^\infty(\mathbb{T}^2)}\right)\|\varphi\|_{H^{s,\gamma}_{g}}\\
&\quad+C_{\gamma, \delta }\left\{1+\|\partial_{xy}U\|^2_{L^\infty(\mathbb{T}^2)}\right\}\|\varphi\|^2_{H^{s,\gamma}_g}+C\|\partial_{xy}^{s}\partial_xP\|_{L^2(\mathbb{T}^2)}
.
\end{aligned}
\end{eqnarray}
\end{Proposition}
$\mathbf{Proof.}$
Applying the operator $\partial_{xy}^{s}$ on $(\ref{1.5})_{1}$ and $(\ref{1.4})_{1}$ respectively, we obtain the following  equations
\begin{equation}\left\{
\begin{array}{ll}
(\partial _{t} +u\partial _{x} +Ku\partial _{y}+w\partial _{z}-\epsilon^2\partial_x^2 -\epsilon^2\partial^2_y -\partial _{z}^{2})\partial_{xy}^{s}\varphi +\partial_{xy}^{s}w\partial _{z}\varphi\\
=- \sum\limits_{0\leq j<s}\binom{s}{j} \partial _{xy}^{s-j}u\partial _{x}\partial _{xy}^{j }\varphi- \sum\limits_{0\leq j<s}\binom{s}{j} \partial _{xy}^{s-j}(Ku) \partial _{y}\partial _{xy}^{j }\varphi\\
\quad- \sum\limits_{1\leq j<s}\binom{s}{j} \partial _{x}^{s-j}w\partial _{z}\partial _{xy}^{j }\varphi + \sum\limits_{0\leq j \leq s}\binom{s}{j} \partial _{x}^{s-j}(\partial_y K u)\partial _{xy}^{j }\varphi   ,\\
(\partial _{t} +u\partial _{x} +Ku\partial _{y}+w\partial _{z}-\epsilon^2\partial _{x}^{2}-\epsilon^2\partial _{y}^{2}-\partial _{z}^{2})\partial_{xy}^{s}(u-U)+\partial_{xy}^{s}w\varphi \\
=- \sum\limits_{0\leq j<s}\binom{s}{j} \partial _{xy}^{s-j}u\partial _{x}\partial _{xy}^{j }(u-U)- \sum\limits_{0\leq j<s}\binom{s}{j} \partial _{xy}^{s-j}(Ku) \partial _{y}\partial _{xy}^{j }(u-U) \\
\quad- \sum\limits_{1\leq j<s}\binom{s}{j} \partial _{xy}^{s-j}w\partial _{xy}^{j }\varphi - \sum\limits_{0\leq j \leq s}\binom{s}{j} \partial _{xy}^{j}(u-U)\partial _{xy}^{s-j }\partial_x U\\
 \quad- \sum\limits_{0\leq j \leq s}\binom{s}{j} \partial _{xy}^{j}(u-U)\partial _{xy}^{s-j }(K\partial_y U ).
\end{array}
 \label{3.1}         \right.\end{equation}
Eliminating $\frac{\partial _{z}\varphi}{\varphi}\times (\ref{3.1})_{2}$ from $ (\ref{3.1})_{1}$, and letting
$a(t,x,y,z)=\frac{\partial _{z}\varphi}{\varphi}$, it follows
\begin{eqnarray}
\begin{aligned}
&
 (\partial _{t} +u\partial _{x} +Ku\partial _{y}+w\partial _{z}-\epsilon^2\partial_x^2-\epsilon^2\partial_y^2-\partial _{z}^{2})g_{s} +\partial_{xy}^{s}(u-U) (\partial _{t} +u\partial _{x} +Ku\partial _{y}+w\partial _{z}-\epsilon^2\partial_x^2-\epsilon^2\partial_y^2-\partial _{z}^{2})a    \\
&= 2\partial_{xy}^{s}\varphi\partial _{z}a+2\epsilon^2 \partial_x\partial^{s}_{xy} (u-U) \partial_x a+2\epsilon^2 \partial_y\partial^{s}_{xy} (u-U) \partial_y a \\
&\quad
- \sum\limits_{j=0}^{s-1}\binom{s}{j} (\partial _{xy}^{s-j}u(\partial _{x}\partial _{xy}^{j }\varphi-a\partial _{x}\partial _{xy}^{j }(u-U))+\partial _{xy}^{s-j}(Ku)(\partial _{y}\partial _{xy}^{j }\varphi-a\partial _{y}\partial _{xy}^{j }(u-U) ))  \\
&\quad - \sum\limits_{j=1}^{s-1}\binom{s}{j}\partial _{xy}^{s-j}w (\partial _{z}\partial _{xy}^{j }\varphi-a \partial _{xy}^{j }\varphi )+\sum\limits_{0\leq j \leq s}\binom{s}{j} \partial _{xy}^{j}(u-U)\partial _{xy}^{s-j }\partial_x U\\
&\quad+a\sum\limits_{0\leq j \leq s}\binom{s}{j} \partial _{xy}^{j}(u-U)\partial _{xy}^{s-j }(K\partial_y U)+ \sum\limits_{0\leq j \leq s}\binom{s}{j} \partial _{x}^{s-j}(\partial_y K u)\partial _{xy}^{j }\varphi . \label{3.2}
\end{aligned}
\end{eqnarray}
Applying the operator $\partial_{z}$ on $(\ref{1.5})_{1}$ yields
\begin{eqnarray}
\begin{aligned}
&(\partial _{t} +u\partial _{x} +Ku\partial _{y} +w\partial _{z})\partial _{z}\varphi \\
&=\epsilon^2(\partial_x^2+\partial_y^2)\partial_z \varphi+\partial _{z}^{3}\varphi-\varphi\partial _{x}\varphi-K\varphi\partial_{y}\varphi+\partial _{x}u\partial _{z}\varphi+\partial _{y}(Ku)\partial _{z}\varphi+\partial_y K(\varphi^2+u\partial_z \varphi).  \label{3.3}
\end{aligned}
\end{eqnarray}
Using $(\ref{3.3})$ and equation $(\ref{1.5})_{1}$, we get
\begin{align}
 &(\partial _{t} +u\partial _{x} +Ku\partial _{y}+w\partial _{z} )a \nonumber \\
 &=\frac{(\partial _{t} +u\partial _{x} +Ku\partial _{y} +w\partial _{z})\partial _{z}\varphi}{\varphi} -\frac{\partial_z \varphi(\partial _{t} +u\partial _{x} +Ku\partial _{y} +w\partial _{z})\varphi}{\varphi^2}\nonumber \\
&=\epsilon^2 \left\{\frac{(\partial_x^2+\partial_y^2)\partial_z \varphi}{\varphi}-a\frac{(\partial_x^2+\partial_y^2)\varphi}{\varphi}\right\}+\left\{\frac{\partial_z^3 \varphi}{\varphi}-a\frac{\partial_z^2 \varphi}{\varphi}\right\}-\partial_x \varphi-K\partial_y \varphi+a\partial_x u+ a\partial_y(Ku)\nonumber \\
&\quad+\partial_y K \varphi,  \label{3.4}
\end{align}
and
\begin{equation}\left\{
\begin{array}{ll}
  \partial _{z}a=\frac{ \partial _{z}^{2}\varphi}{\varphi} -\frac{\partial _{z}w \partial _{z}\varphi}{\varphi^{2}},   \\
 \partial _{z}^{2}a=\frac{\partial _{z}^{3}\varphi}{\varphi}-a\frac{\partial _{z}^{2}\varphi}{\varphi}- 2a \partial _{z}a, \\
 \partial _{x}^{2}a=\frac{\partial _{x}^{2}\partial_z\varphi}{\varphi}-a\frac{\partial _{x}^{2}\varphi}{\varphi}- 2\frac{\partial_x \varphi}{\varphi} \partial _{x}a,\\
 \partial _{y}^{2}a=\frac{\partial _{y}^{2}\partial_z\varphi}{\varphi}-a\frac{\partial _{y}^{2}\varphi}{\varphi}- 2\frac{\partial_y \varphi}{\varphi} \partial _{y}a.
\end{array}
 \label{3.5}         \right.\end{equation}
 Substituting (\ref{3.5}) into (\ref{3.4}), we obtain an equation for $a$:
 \begin{eqnarray}
\begin{aligned}
&(\partial _{t} +u\partial _{x} +Ku\partial _{y}+w\partial _{z}-\epsilon^2\partial_x^2-\epsilon^2\partial_y^2-\partial _{z}^{2} )a\\
&=2\epsilon^2(\frac{\partial_x \varphi}{\varphi} \partial _{x}a+\frac{\partial_y \varphi}{\varphi} \partial _{y}a)+2a \partial_z a-g_{1x}+a\partial_x U -Kg_{1y}+K a\partial_y U\\
&\quad+a\partial_y K u+\partial_y K \varphi.
\label{3.444}
\end{aligned}
\end{eqnarray}

Now inserting $(\ref{3.444})$ and $(\ref{3.5})_{2}$ into $(\ref{3.2})$ yields
\begin{eqnarray}
\begin{aligned}
&
 (\partial _{t} +u\partial _{x} +Ku\partial _{y}+w\partial _{z}-\epsilon^2\partial_x^2-\epsilon^2\partial_y^2-\partial _{z}^{2})g_{s}    \\
&= 2g_{s}\partial _{z}a+2\epsilon^2\partial_x a\left( \partial_x\partial^{s}_{xy} (u-U)-\frac{\partial_x \varphi}{\varphi}\partial^{s}_{xy} (u-U)\right) +2\epsilon^2\partial_y a\left( \partial_y\partial^{s}_{xy} (u-U)-\frac{\partial_y \varphi}{\varphi}\partial^{s}_{xy} (u-U)\right) \\
&\quad
- \sum\limits_{j=1}^{s-1}\binom{s}{j} g_{j+1}\partial _{xy}^{s-j}u-g_{x1}\partial_{xy}^s U- \sum\limits_{j=1}^{s-1}\binom{s}{j} g_{j+1}\partial _{xy}^{s-j}(Ku)-  g_{1y}\partial _{xy}^{s}(Ku)+g_{1y}(K\partial_{xy}^{s}(u-U))   \\
&\quad - \sum\limits_{j=1}^{s-1}\binom{s}{j}\partial _{xy}^{s-j}w (\partial _{z}\partial _{xy}^{j }\varphi-a \partial _{xy}^{j }\varphi )+a\sum\limits_{0\leq j \leq s-1}\binom{s}{j} \partial _{xy}^{j}(u-U)\partial _{xy}^{s-j }\partial_x U\\
&\quad+a\sum\limits_{0\leq j \leq s-1}\binom{s}{j} \partial _{xy}^{j}(u-U)\partial _{xy}^{s-j }(K\partial_y U)+ \sum\limits_{0\leq j \leq s}\binom{s}{j} \partial _{x}^{s-j}(\partial_y K u)\partial _{xy}^{j }\varphi \\
&\quad -( a\partial_y K u+\partial_y K \varphi)\partial_{xy}^{s}(u-U). \label{3.6}
\end{aligned}
\end{eqnarray}
 Multiplying $(\ref{3.6})$ by $(1+z)^{2\gamma}g_{s}$, and integrating it over $\mathbb{T}^{2}\times \mathbb{R}^{+}$, we arrive at
\begin{align}
&\frac{1}{2}\frac{d}{dt}\|(1+z)^{ \gamma}g_{s}\|_{L^{2}}^{2}  +\|(1+z)^{\gamma}\partial _{z}g_{s}\|_{L^{2}} ^{2}+\epsilon^2\|(1+z)^{\gamma}\partial _{x}g_{s}\|_{L^{2}} ^{2}+\epsilon^2\|(1+z)^{\gamma}\partial _{y}g_{s}\|_{L^{2}} ^{2} \nonumber \\
&=\iint_{\mathbb{T}^{2} } g_{s}\partial _{z}g_{s}dxdy | _{z=0}-2\gamma \iiint(1+z)^{2\gamma-1}g_{s}\partial _{z}  g_{s}
+\gamma \iiint(1+z)^{2\gamma-1} w  |g_{s}|^{2}  \nonumber\\
&\quad +2\iiint(1+z)^{2\gamma}|g_{s}|^{2}\partial _{z}a+2\epsilon^2\iiint(1+z)^{2\gamma}g_{s}\partial_x a\left( \partial_x\partial^{s}_{xy} (u-U)-\frac{\partial_x \varphi}{\varphi}\partial^{s}_{xy} (u-U)\right)\nonumber\\
&\quad +2\epsilon^2\iiint(1+z)^{2\gamma}g_s\partial_y a\left( \partial_y\partial^{s}_{xy} (u-U)-\frac{\partial_y \varphi}{\varphi}\partial^{s}_{xy} (u-U)\right) \nonumber\\
&\quad - \sum\limits_{j=1}^{s-1}\binom{s}{j} \iiint(1+z)^{2\gamma}g_s g_{j+1}\partial _{xy}^{s-j}u-\iiint(1+z)^{2\gamma}g_s g_{x1}\partial_{xy}^s U\nonumber\\
&\quad- \sum\limits_{j=1}^{s-1}\binom{s}{j} \iiint(1+z)^{2\gamma}g_sg_{j+1}\partial _{xy}^{s-j}(Ku)- \iiint(1+z)^{2\gamma}g_s g_{1y}\partial _{xy}^{s}(Ku)\nonumber\\
&\quad+\iiint(1+z)^{2\gamma}g_s g_{1y}(K\partial_{xy}^{s}(u-U))- \sum\limits_{j=1}^{s-1}\binom{s}{j}\iiint(1+z)^{2\gamma}g_s\partial _{xy}^{s-j}w (\partial _{z}\partial _{xy}^{j }\varphi-a \partial _{xy}^{j }\varphi )\nonumber\\
&\quad +\sum\limits_{0\leq j \leq s-1}\binom{s}{j} \iiint(1+z)^{2\gamma}g_s a\partial _{xy}^{j}(u-U)\partial _{xy}^{s-j }\partial_x U\nonumber\\
&\quad+\sum\limits_{0\leq j \leq s-1}\binom{s}{j}\iiint(1+z)^{2\gamma}g_s a \partial _{xy}^{j}(u-U)\partial _{xy}^{s-j }(K\partial_y U)\nonumber\\
&\quad+ \sum\limits_{0\leq j \leq s}\binom{s}{j}\iiint(1+z)^{2\gamma}g_s \partial _{x}^{s-j}(\partial_y K u)\partial _{xy}^{j }\varphi -\iiint(1+z)^{2\gamma}g_s( a\partial_y K u+\partial_y K \varphi)\partial_{xy}^{s}(u-U)
, \label{3.7}
\end{align}
which can be obtained by integration and by using the boundary condition, i.e.,
\begin{eqnarray}
\begin{aligned}
& \iiint(1+z)^{2\gamma}g_{s}(u\partial _{x} +Ku\partial _{y}+w\partial _{z} )g_{s} \\
&=\frac{1}{2} \iiint(1+z)^{2\gamma}\left\{u(|g_{s}|^{2})_{x}+ u(|g_{s}|^{2}) _{y}+w (|g_{s}|^{2}) _{z}\right\}
 =-\gamma \iiint(1+z)^{2\gamma-1}w |g_{s}|^{2},
 \label{3.8}
\end{aligned}
\end{eqnarray}
and
\begin{align}
 \iiint(1+z)^{2\gamma}g_{s}\partial _{z}^{2} g_{s}   &=-\|(1+z)^{\gamma}\partial _{z}g_{s}\|_{L^{2}}-\iint_{\mathbb{T}^{2}} g_{s}\partial _{z}g_{s}dxdy | _{z=0}
 -2\gamma \iiint(1+z)^{2\gamma-1}g_{s}\partial _{z}  g_{s}.
 \label{3.9}
\end{align}
By definition of $H^{s,\gamma}_{\sigma,\delta}$, we have
\begin{equation}\left\{
\begin{array}{ll}
|\varphi|^{-1}\leq \delta^{-1}(1+z)^{\sigma},   \\
| \partial _{x} \varphi|, |\partial _{y} \varphi| \leq \delta^{-1}(1+z)^{-\sigma },  \\
|\partial _{z}\varphi|, |\partial _{x} \partial _{z}\varphi|, |\partial _{y} \partial _{z}\varphi| \leq \delta^{-1}(1+z)^{-\sigma-1},\\
|\partial _{z}^{2}\varphi|  \leq \delta^{-1}(1+z)^{-\sigma-2}.
\end{array}
 \label{2.60}         \right.\end{equation}
We estimate equation $(\ref{3.7})$ by terms. Firstly,  we have the fact that
\begin{align}
\partial _{z}g_{s} \big{|} _{z=0}&=\left(\partial_{xy}^{s}\partial _{z}\varphi-\frac{\partial _{z}\varphi}{\varphi} \partial_{xy}^{s}\varphi -\frac{\partial^{2} _{z}\varphi}{\varphi} \partial_{xy}^{s}(u-U) +\frac{\partial_{z}\varphi\partial_{z}\varphi}{\varphi^{2}} \partial_{xy}^{s}(u-U) \right) \big{|} _{z=0} \nonumber\\
&=\partial_{xy}^{s}\partial _{x}P-ag_{s} \big{|} _{z=0} +\frac{\partial^{2} _{z}\varphi}{\varphi} \partial_{xy}^{s}U \big{|} _{z=0}.  \label{3.10}
\end{align}
Then according to the trace theorem, together with the facts that $\|a\|_{L^{\infty}}\leq \delta^{-2}$ and $\|\partial _{z}a\|_{L^{\infty}}\leq \delta^{-2}+\delta^{-4}$, we get
\begin{eqnarray}
\left|\iint_{\mathbb{T}^{2} } g_{s}\partial _{z}g_{s}dxdy | _{z=0}  \right|\leq \frac{1}{4}\|(1+z)^{\gamma}\partial_{z} g_{s}\|_{L^{2}}^{2}+C_{\delta}\left\{1+\|\partial_{xy}U\|^2_{L^\infty(\mathbb{T}^2)}\right\}\|\varphi\|^2_{H^{s,\gamma}_g}+C\|\partial_{xy}^{s}\partial_xP\|_{L^2(\mathbb{T}^2)},   \label{3.11}
\end{eqnarray}
Using now the  H$\ddot{o}$lder inequality,  we obtain
\begin{eqnarray}
\left| 2\gamma \iiint(1+z)^{2\gamma-1}g_{s}\partial _{z}  g_{s}  \right|
 &\leq& 2\gamma \|\frac{1}{1+z }  \|_{L^{\infty}} \|(1+z)^{\gamma} g_{s}\|_{L^{2}} \|(1+z)^{\gamma}\partial _{z}g_{s}\|_{L^{2}} \nonumber \\
& \leq & C_{ \gamma } \|\varphi\|^2_{H^{s,\gamma}_g}+\frac{1}{4}\|(1+z)^{\gamma}\partial_{z} g_{s}\|_{L^{2}}^{2} ,  \label{3.12}
\end{eqnarray}
and using Lemma \ref{y4.4},
\begin{align}
\left|   \gamma \iiint(1+z)^{2\gamma-1} w  |g_{s}|^{2} \right| & \leq \|\frac{w}{1+z}\|_{L^{\infty}}\|(1+z)^{ \gamma}g_{s}\|_{L^{2}}^{2}\nonumber\\
&\leq  C_{s, \gamma,\sigma,\delta  }\left(1+\|\varphi \|_{H^{s,\gamma}_{g}}+\|\partial_{xy}U\|^2_{L^\infty(\mathbb{T}^2)}\right)\|\varphi \|_{H^{s,\gamma}_{g}} ^2.  \label{3.13}
\end{align}

Thus, it follows from to $(\ref{2.60})$ and $(\ref{3.5})$ that
\begin{eqnarray}
\left| 2\iiint(1+z)^{2\gamma}|g_{s}|^{2}\partial _{z}a  \right|
  \leq  2  \|\partial _{z}a  \|_{L^{\infty}} \|(1+z)^{\gamma} g_{s}\|_{L^{2}} ^{2}\leq C_{ \delta }  \|\varphi \|_{H^{s,\gamma}_{g}} ^2,
\label{3.15}
\end{eqnarray}
and by Lemma \ref{y4.2},
\begin{eqnarray}
&&\left| 2\epsilon^2\iiint(1+z)^{2\gamma}g_{s} \partial_x a\left( \partial_x\partial^{s}_{xy} (u-U)-\frac{\partial_x \varphi}{\varphi}\partial^{s}_{xy} (u-U)\right)   \right| \nonumber \\
&&\leq  2\epsilon^2C_{\delta} \|(1+z)^{\gamma} g_{s}\|_{L^{2}}\left(\|(1+z)^{\gamma-1}\partial_{xy}^s\partial_x(u-U)\|_{L^2}+\|(1+z)^{\gamma-1}\partial_{x}^s(u-U)\|_{L^2}\right),
\end{eqnarray}
which, together with $(\ref{2.60})$, implies the fact that $\|(1+z)(\partial _{x}a+\partial _{y}a)\|_{L^{\infty}} \leq \delta^{-2}+\delta^{-4}$.
In addition, by using Lemma \ref{y4.1} and  the following relation
\begin{eqnarray}
\varphi\partial_z\left(\frac{\partial_x\partial_{xy}^{s}(u-U)}{\varphi}\right)=g_{(s+1)x}=\partial_{x} g_s +\partial_x a\partial_{xy}^s(u-U),
\end{eqnarray}
we get
\begin{align}
&\|(1+z)^{\gamma-1}\ \partial_x\partial^{s}_{xy}(u-U)\|_{L^2} \nonumber\\
&\leq  \delta^{-1}\|(1+y)^{\gamma-\sigma-1}\frac{\partial_x\partial_{xy}^{s}(u-U)}{\varphi}\|_{L^2}\nonumber\\
&\leq C_{\gamma, \sigma, \delta}\left(\|\partial_{x}\partial_{xy}^{s}U\|_{L^{2}(\mathbb{T})}+\left\|(1+z)^{\gamma}\varphi\partial_z\left(\frac{\partial_x\partial_{xy}^{s}(u-U)}{\varphi}\right)\right\|_{L^2}\right)\nonumber\\
&\leq C_{\gamma, \sigma, \delta}\left(\|\partial_{x}\partial_{xy}^{s}U\|_{L^{2}(\mathbb{T}^2)}+\|(1+z)^{\gamma} \partial_x g_{s}\|_{L^2}+\|(1+z)\partial_x a\|_{L^{\infty}}\|(1+z)^{\gamma-1}\partial_{xy}^s(u-U)\|_{L^2}\right) .
\end{align}
Then, we can obtain
\begin{eqnarray}
&&\left| 2\epsilon^2\iiint(1+z)^{2\gamma}g_{s} \partial_x a\left( \partial_x\partial^{s}_{xy} (u-U)-\frac{\partial_x \varphi}{\varphi}\partial^{s}_{xy} (u-U)\right)   \right| \nonumber \\
&&\leq  \epsilon^2C_{\gamma, \sigma, \delta} \left(\|\varphi\|_{H^{s,\gamma}_{g}}+\|\partial_{xy}^s U\|_{L^2(\mathbb{T}^2)}+\|\partial_{x}\partial_{xy}^{s}U\|_{L^{2}(\mathbb{T}^2)}+\|(1+z)^{\gamma} \partial_x g_{s}\|_{L^2}\right)\|\varphi\|_{H^{s,\gamma}_{g}}.
\label{3.1555}
\end{eqnarray}
Similarly, we also have
\begin{eqnarray}
&&\left| 2\epsilon^2\iiint(1+z)^{2\gamma}g_{s} \partial_y a\left( \partial_y\partial^{s}_{xy} (u-U)-\frac{\partial_y \varphi}{\varphi}\partial^{s}_{xy} (u-U)\right)   \right| \nonumber \\
&&\leq  \epsilon^2C_{\gamma, \sigma, \delta} \left(\|\varphi\|_{H^{s,\gamma}_{g}}+\|\partial_{xy}^s U\|_{L^2(\mathbb{T}^2)}+\|\partial_{y}\partial_{xy}^{s}U\|_{L^{2}(\mathbb{T}^2)}+\|(1+z)^{\gamma} \partial_y g_{s}\|_{L^2}\right)\|\varphi\|_{H^{s,\gamma}_{g}}.
\label{3.15555}
\end{eqnarray}
Now we estimate some of the remaining complex terms. Thus we need to discuss the classification to complete the estimation. For $i \leq s-3$, we use the often used inequality to derive
\begin{align}
\|(1+z)^{\gamma}g_i\|_{L^\infty} &\leq \|(1+z)^{\gamma}\partial^i_{xy}\varphi\|_{L^\infty}+\|(1+z)^{\gamma-1}a\|_{L^\infty}\|\partial_{xy}^i (u-U)\|_{L^\infty} \nonumber \\
&\leq C_{\delta} \left(\|(1+z)^{\gamma+\alpha_3}\partial^i_{xy}\varphi\|_{L^\infty}+\|\partial_{xy}^i u\|_{L^\infty}+\|\partial_{xy}^i U\|_{L^\infty}\right)\nonumber \\
&\leq C_{s,\gamma, \sigma, \delta}\left(\|\varphi\|_{H^{s,\gamma}_{g}}+\|\partial_{xy}^s U\|_{L^2(\mathbb{T}^2)}\right).
\end{align}
By Lemma \ref{y4.4}, we obtain
\begin{align}
\left| \iiint(1+z)^{2\gamma}g_s g_{j+1}\partial _{xy}^{s-j}u\right|
&\leq \left\{
\begin{aligned}
&\|(1+z)^{\gamma}  g_{s}\|_{L^2}\|(1+z)^{\gamma}  g_{j+1}\|_{L^2}\|\partial _{xy}^{s-j}u\|_{L^\infty}, ~j=[2, s-1],\\
&\|(1+z)^{\gamma}  g_{s}\|_{L^2}\|(1+z)^{\gamma}  g_2\|_{L^\infty}\|\partial _{xy}^{s-1}(u-U)+\partial _{xy}^{s-1}U\|_{L^2}, ~j=1.\\
\end{aligned}
\right.\nonumber \\
&\leq C_{s,\gamma, \sigma, \delta}\left(\|\varphi\|_{H^{s,\gamma}_{g}}+\|\partial_{xy}^s U\|_{L^2(\mathbb{T}^2)}\right)^2\|\varphi\|_{H^{s,\gamma}_{g}},
\label{3.155}
\end{align}
for $j=1,2,\cdots, s-1$, and
\begin{eqnarray}
\left|\iiint(1+z)^{2\gamma}g_s g_{x1}\partial_{xy}^s U\right| \leq C_{s,\gamma, \sigma, \delta} \left(\|\varphi\|_{H^{s,\gamma}_{g}}+\|\partial_{xy}^s U\|_{L^2(\mathbb{T}^2)}\right)\|\partial_{xy}^s U\|_{L^2(\mathbb{T}^\infty)}\|\varphi\|_{H^{s,\gamma}_{g}}.
\end{eqnarray}
Similarly,
\begin{eqnarray}
\begin{aligned}
&\left| \sum\limits_{j=1}^{s-1}\binom{s}{j} \iiint(1+z)^{2\gamma}g_sg_{j+1}\partial _{xy}^{s-j}(Ku)\right|\leq C_{s,\gamma, \sigma, \delta}\|\partial_{xy}^{s-1} K\|_{L^\infty}\left(\|\varphi\|_{H^{s,\gamma}_{g}}+\|\partial_{xy}^s U\|_{L^2(\mathbb{T}^2)}\right)^2\|\varphi\|_{H^{s,\gamma}_{g}}
\end{aligned}
\end{eqnarray}
and
\begin{align}
&\left|- \iiint(1+z)^{2\gamma}g_s g_{1y}\partial _{xy}^{s}(Ku)+\iiint(1+z)^{2\gamma}g_s g_{1y}(K\partial_{xy}^{s}(u-U))\right|\nonumber\\
&\leq C_{s,\gamma, \sigma, \delta}\|\partial_{xy}^{s-1} K\|_{L^\infty}\|(1+z)^{\gamma}  g_{s}\|_{L^2}\|(1+z)^{\gamma}  g_{1y}\|_{L^\infty}\|-\partial _{xy}^{s}(Ku)+(K\partial_{xy}^{s}(u-U))\|_{L^2}\nonumber\\
& \leq C_{s,\gamma, \sigma, \delta}\|\partial_{xy}^{s-1} K\|_{L^\infty}\left(\|\varphi\|_{H^{s,\gamma}_{g}}+\|\partial_{xy}^s U\|_{L^2(\mathbb{T}^2)}\right)^2\|\varphi\|_{H^{s,\gamma}_{g}}.
\end{align}
It should be noted that  the highest order of $u$ will not reach $s$, due to cancelling each other out,  in the above  inequality.\\
Last, for $j=3,\cdots, s-1$,
\begin{eqnarray}
\left|   \iiint(1+z)^{2\gamma}g_s\partial _{xy}^{s-j}w (\partial _{z}\partial _{xy}^{j }\varphi-a \partial _{xy}^{j }\varphi )  \right| \leq C_{s,\gamma, \sigma, \delta}\left(\|\varphi\|_{H^{s,\gamma}_{g}}+\|\partial_{xy}^s U\|_{L^2(\mathbb{T}^2)}\right)\|\varphi\|_{H^{s,\gamma}_{g}}^2
\label{3.17}
\end{eqnarray}
and for $j=1, 2$,
\begin{eqnarray}
\left|   \iiint(1+z)^{2\gamma}g_s\partial _{xy}^{s-j}w (\partial _{z}\partial _{xy}^{j}\varphi-a \partial _{xy}^{j}\varphi )  \right| \leq C_{s,\gamma, \sigma, \delta}\left(\|\varphi\|_{H^{s,\gamma}_{g}}+\|\partial_{xy}^s U\|_{L^2(\mathbb{T}^2)}\right)\|\varphi\|_{H^{s,\gamma}_{g}}^2.
\label{3.18}
\end{eqnarray}
Similar to estimates used in \cite{3} on the 2D Prandtl equations,  we directly here give the above two estimates.
For $j=1,2,\cdots, s-1$, by Lemma \ref{y4.4} and using the fact  $\|(1+z) a \|_{L^{\infty}} \leq \delta^{-2}$, it follows
\begin{eqnarray}
&& \left|    \iiint(1+z)^{2\gamma}g_s a\partial _{xy}^{j}(u-U)\partial _{xy}^{s-j }\partial_x U  \right|
 \nonumber \\
&& \leq    C \|(1+z)^{\gamma}g_s\|_{L^2}\|(1+z)a\|_{L^\infty}\|(1+z)^{\gamma-1}\partial_{xy}^{j}(u-U)\|_{L^2}\|\partial _{x}\partial _{xy}^{s-j}U\|_{L^\infty(\mathbb{T}^2)}  \nonumber \\
&& \leq C_{s,\gamma, \sigma, \delta}\|\partial _{x}\partial _{xy}^{s}U\|_{L^\infty(\mathbb{T}^2)} \left(\|\varphi\|_{H^{s,\gamma}_{g}}+\|\partial_{xy}^s U\|_{L^2(\mathbb{T}^2)}\right)\|\varphi\|_{H^{s,\gamma}_{g}},
\label{3.19}
\end{eqnarray}
and
\begin{eqnarray}
&& \left|    \iiint(1+z)^{2\gamma}g_s a \partial _{xy}^{j}(u-U)\partial _{xy}^{s-j }(K\partial_y U)  \right|
 \nonumber \\
&& \leq C_{s,\gamma, \sigma, \delta}\|\partial _{xy}^{s}(K\partial_yU)\|_{L^\infty(\mathbb{T}^2)} \left(\|\varphi\|_{H^{s,\gamma}_{g}}+\|\partial_{xy}^s U\|_{L^2(\mathbb{T}^2)}\right)\|\varphi\|_{H^{s,\gamma}_{g}}.
\label{3.199}
\end{eqnarray}
For $j=1,2,\cdots, s$, we achieve
\begin{align}
\left| \iiint(1+z)^{2\gamma}g_s \partial _{x}^{s-j}(\partial_y K u)\partial _{xy}^{j }\varphi\right|
&\leq C_{\delta}\|\varphi\|_{H^{s,\gamma}_{g}}\left\{
\begin{aligned}
&\|(1+z)^{\gamma} \partial _{xy}^{j }\varphi\|_{L^2}\|\partial _{x}^{s-j}(\partial_y K u)\|_{L^\infty}, ~j \geq 2,\\
&\|(1+z)^{\gamma} \partial _{xy}^{j }\varphi\|_{L^\infty}\|\partial _{x}^{s-j}(\partial_y K u)\|_{L^2}, ~j=0,1.\\
\end{aligned}
\right.\nonumber \\
&\leq C_{s,\gamma, \sigma, \delta}\|\partial _{xy}^{s-2}K\|_{L^\infty(\mathbb{T}^2)}\left(\|\varphi\|_{H^{s,\gamma}_{g}}+\|\partial_{xy}^s U\|_{L^2(\mathbb{T}^2)}\right)^2\|\varphi\|_{H^{s,\gamma}_{g}},\label{3.1999}
\end{align}
and
\begin{align}
-\iiint(1+z)^{2\gamma}g_s( a\partial_y K u+\partial_y K \varphi)\partial_{xy}^{s}(u-U) &\leq  C\|\partial_y K \|_{L^\infty}\|(1+z)^{\gamma}g_s\|_{L^2}\|(1+z)^{\gamma-1}\partial_{xy}^{s}(u-U)\|_{L^2}
\nonumber \\
&\quad \times (\|(1+z)a\|_{L^\infty}\|u\|_{L^\infty}+\|(1+z)^{1-\gamma}\|_{L^\infty}\|(1+z)^{\gamma}\varphi\|_{L^\infty})\nonumber\\
&\leq C_{s,\gamma, \sigma, \delta} \|\partial_y K \|_{L^\infty}\left(\|\varphi\|_{H^{s,\gamma}_{g}}+\|\partial_{xy}^s U\|_{L^2(\mathbb{T}^2)}\right)^2\|\varphi\|_{H^{s,\gamma}_{g}}.
\label{3.20}
\end{align}

Hence combining (\ref{3.11})-(\ref{3.15}), (\ref{3.1555})-(\ref{3.15555}), (\ref{3.155})-(\ref{3.20}) with (\ref{3.7}) leads to
\begin{eqnarray}
\begin{aligned}
\frac{d}{dt}\|(1+z)^{ \gamma}g_{s}\|_{L^{2}}^{2}& \leq  C_{s,\gamma, \sigma, \delta}\left(1+ \|\partial_{xy}^s K \|_{L^\infty}\right)\left(1+\|\varphi\|_{H^{s,\gamma}_{g}}+\|\partial_{xy}^s U\|_{L^\infty(\mathbb{T}^2)}\right)\\
&\quad \times\left(\|\varphi\|_{H^{s,\gamma}_{g}}+\|\partial_{xy}^{s+1} U\|_{L^\infty(\mathbb{T}^2)}\right)\|\varphi\|_{H^{s,\gamma}_{g}}\\
&\quad+C_{\gamma, \delta }\left\{1+\|\partial_{xy}U\|^2_{L^\infty(\mathbb{T}^2)}\right\}\|\varphi\|^2_{H^{s,\gamma}_g}+C\|\partial_{xy}^{s}\partial_xP\|_{L^2(\mathbb{T}^2)}
.
\end{aligned}
\end{eqnarray}
$\hfill$$\Box$

\subsection{Weighted $H^s$ estimate on $\varphi$ }
Now we can derive the weighted $H^s$ estimate on $\varphi$ by employing Propositions \ref{p2.1}-\ref{p2.2}.
\begin{Proposition}\label{p2.3}
Under the same assumption of Proposition \ref{p2.1}, we have the following estimate:
\begin{eqnarray}
\begin{aligned}
\|\varphi\|_{H^{s,\gamma}_{g}}^2 &\leq \left\{\|\varphi_{0}\|_{H^{s,\gamma}_{g}}^2+\int_0^t F(\tau)d \tau\right\}\\
&\quad \times \left\{1-C(s-1)\left(\|\varphi_{0}\|_{H^{s,\gamma}_{g}}^2+\int_0^t F(\tau) d\tau\right)^{s-1}t\right\}^{-\frac{1}{s-1}},
\label{3.2222}
\end{aligned}
\end{eqnarray}
where $C>0$ is a constant independent of $\epsilon$ and $t$. The function $F(t)$ is expressed by
\begin{eqnarray}
\begin{aligned}
F(t)=\mathcal{P}(\|\partial_{xy}^{s+1}U\|_{{L^\infty}(\mathbb{T}^2)},
\|\partial_{xy}^{s}K\|_{{L^\infty}(\mathbb{T}^2)})+C \sum\limits_{l=0}^{\frac{s}{2}}\|\partial^{l}_t\partial_xP\|_{H^{s-2l}(\mathbb{T}^2)}^{2},
\end{aligned}
\end{eqnarray}
and $\mathcal{P}(\cdot)$ denotes a polynomial.
\end{Proposition}
$\mathbf{Proof.}$
According to Propositions \ref{p2.1}-\ref{p2.2}, we deduce from the definition of $\|\cdot\|_{H^{s,\gamma}_g}$ that
\begin{align}
\frac{d}{dt} \|\varphi\|_{H^{s,\gamma}_{g}}^2 &\leq C_{s,\gamma, \sigma, \delta} \|\varphi\|_{H^{s,\gamma}_{g}}^{2s}+\sum\limits_{l=0}^{\frac{s}{2}}\|\partial^{l}_t\partial_xP\|_{H^{s-2l}(\mathbb{T}^2)}^{2}+\mathcal{P}(\|\partial_{xy}^{s+1}U\|_{{L^\infty}(\mathbb{T}^2)},
\|\partial_{xy}^{s}K\|_{{L^\infty}(\mathbb{T}^2)})
\end{align}
and hence, it follows by using the comparison principle of ordinary differential equations that
\begin{eqnarray}
\begin{aligned}
\|\varphi\|_{H^{s,\gamma}_{g}}^2 &\leq \left\{\|\varphi_0\|_{H^{s,\gamma}_{g}}^2+\int_0^t F(\tau)d \tau\right\}\\
&\quad \times \left\{1-C(s-1)\left(\|\varphi_0\|_{H^{s,\gamma}_{g}}^2+\int_0^t F(\tau) d\tau\right)^{s-1}t\right\}^{-\frac{1}{s-1}},
\end{aligned}
\end{eqnarray}
provided that
\begin{eqnarray}
\begin{aligned}
1-C(s-1)\left(\|\varphi_0\|_{H^{s,\gamma}_{g}}^2+\int_0^t F(\tau) d\tau\right)^{s-1}t > 0.
\end{aligned}
\end{eqnarray}
This proves Proposition \ref{p2.3}.
$\hfill$$\Box$

\subsection{Weighted $L^{\infty}$ estimates on lower order terms  }
In this subsection, we will estimate the weighted $L^{\infty}$ on $D^{\alpha} \varphi$ for $|\alpha|\leq 2$ by using the classical maximum principles. More precisely, we will derive two parts: an $L^\infty$-estimate on  $I:=\sum \limits_{|\alpha|\leq 2} \left|(1+z)^{\sigma+\alpha_{3}}D^{\alpha}\varphi\right|^2$ and a lower bound estimate on $B_{(0,0,0)}:=(1+z)^{\sigma} \varphi$. Here,
\begin{eqnarray*}
\begin{aligned}
I:=\sum \limits_{|\alpha|\leq 2} \left|(1+z)^{\sigma+\alpha_{3}}D^{\alpha}\varphi\right|^2,\;B_{\alpha}:=(1+z)^{\sigma+\alpha_{3}}D^{\alpha}\varphi,
\end{aligned}
\end{eqnarray*}
and $ B_{(0,0, 0)}=(1+z)^{\sigma}\varphi.$

\begin{Lemma}\label{y2.3}
Under the same assumption of Proposition \ref{p2.1}, we have the following estimates: \\
For any $s \geq 5$,
\begin{eqnarray}
\begin{aligned}
\|I(t)\|_{L^\infty(\mathbb{T}^2\times\mathbb{R}^+)}\leq \max \left\{           \|I(0)\|_{L^\infty(\mathbb{T}^2\times\mathbb{R}^+)},6C^2W(t)^2
\right\} e^{C\left(1+G_1(t)\right)t},
\label{33311}
\end{aligned}
\end{eqnarray}
and for any $s\geq 7$,
\begin{eqnarray}
\begin{aligned}
\|I(t)\|_{L^\infty(\mathbb{T}^2\times\mathbb{R}^+)}\leq  \left\{           \|I(0)\|_{L^\infty(\mathbb{T}^2\times\mathbb{R}^+)}+C[1+Y(t)W(t)]W(t)^2t
\right\} e^{C\left(1+G_1(t)\right)t}.
\label{3.32}
\end{aligned}
\end{eqnarray}
In addition, if $s \geq 5$, we also have,
\begin{eqnarray}
\begin{aligned}
\mathop{\min}\limits_{\mathbb{T}^2\times\mathbb{R}^+}(1+z)^{\sigma}\varphi(t)\geq \max \left\{         1-C\left(1+G_2(t)\right)te^{C\left(1+G_2(t)\right)t}
\right\}\cdot\min \left\{\mathop{\min}\limits_{\mathbb{T}^2\times\mathbb{R}^+}(1+z)^{\sigma}\varphi_0-CW(t)t
\right\},
\label{3.33}
\end{aligned}
\end{eqnarray}
where constant $C>0$ depends on $s, \gamma,\sigma$, and $\delta$ only. The functions $G_1$, $G_2$, W and Y : $[0,T] \rightarrow \mathbb{R}^+$ are respectively defined by
\begin{eqnarray}
\begin{aligned}
G_1(t):&=\mathop{\sup}\limits_{ [0,t]}\|\partial_{xy}^{2}K(\tau)\|_{{L^{\infty}}(\mathbb{T}^2)}+\mathop{\sup}\limits_{ [0,t]}\|w(\tau)\|_{H^{s,\gamma}_{g}}+\mathop{\sup}\limits_{ [0,t]}\|\partial_{xy}^{s}U(\tau)\|_{{L^2}(\mathbb{T}^2)}\\
&\quad+\mathop{\sup}\limits_{ [0,t]}\|\partial_{xy}^{2}K(\tau)\|_{{L^{\infty}}(\mathbb{T}^2)}\left(\mathop{\sup}\limits_{ [0,t]}\|w(\tau)\|_{H^{s,\gamma}_{g}}+\mathop{\sup}\limits_{ [0,t]}\|\partial_{xy}^{s}U(\tau)\|_{{L^2}(\mathbb{T}^2)}\right) ,\\
G_2(t):&=\mathop{\sup}\limits_{ [0,t]}\|w(\tau)\|_{H^{s,\gamma}_{g}}+\mathop{\sup}\limits_{ [0,t]}\|\partial_{xy}^{s}U(\tau)\|_{{L^2}(\mathbb{T}^2)},
\label{2.012}
\end{aligned}
\end{eqnarray}
and
\begin{eqnarray}
  \ W(t):=\mathop{\sup}\limits_{ [0,t]}\|w(\tau)\|_{H^{s,\gamma}_{g}}, \ \ \ \ Y(t):=\mathop{\sup}\limits_{ [0,t]}\|\partial_{xy}K(\tau)\|_{{L^{\infty}}(\mathbb{T}^2)}.
\label{2.0122}
\end{eqnarray}
\end{Lemma}
$\mathbf{Proof.}$
By direct computations, we obtain
\begin{eqnarray}
\begin{aligned}
\left(\partial_t+u\partial_x+Ku\partial_y+w\partial_z-\partial^2_z-\epsilon^2\partial^2_x-\epsilon^2\partial^2_y-\partial_y K u\right)B_{\alpha}=\sum\limits_{i=1}^{3} S_{i},
\label{3.34}
\end{aligned}
\end{eqnarray}
where
\begin{eqnarray*}
\begin{aligned}
S_1=\left(\frac{\sigma+\alpha_3}{1+z}w+\frac{(\sigma+\alpha_3)(\sigma+\alpha_3-1)}{(1+z)^2}\right)B_{\alpha},\ \ S_2=-\frac{2(\sigma+\alpha_3)}{1+z}\partial_z B_{\alpha},
\end{aligned}
\end{eqnarray*}
and
\begin{eqnarray*}
\begin{aligned}
S_3=-\sum\limits_{0< \beta \leq \alpha}\binom{\alpha}{\beta}\left\{(1+y)^{\beta_3}\left(D^{\beta}u B_{\alpha-\beta+e_1}+D^{\beta}(Ku) B_{\alpha-\beta+e_2}+\frac{D^{\beta} w B_{\alpha-\beta+e_3}}{1+z}+D^{\beta}(\partial_{y} Ku) B_{\alpha-\beta}\right)\right\}.
\end{aligned}
\end{eqnarray*}
Multiplying the equation (\ref{3.34}) by $2B_{\alpha}$  gives
\begin{eqnarray}
&&\left(\partial_t+u\partial_x+Ku\partial_y+w\partial_z-\partial^2_z-\epsilon^2\partial^2_x-\epsilon^2\partial^2_y-2\partial_y K u\right)I\nonumber\\
&&=2\sum \limits_{|\alpha|\leq 2} \left(|S_1 B_{\alpha}|+|S_2 B_{\alpha}|+|S_3 B_{\alpha}|-\epsilon^2|\partial_x B_{\alpha}|^2-\epsilon^2|\partial_y B_{\alpha}|^2-|\partial_z B_{\alpha}|^2\right)\nonumber\\
&&\leq  C_{s, \gamma,\sigma,\delta  }(1+\|\partial_{xy}^{2}K\|_{{L^{\infty}}(\mathbb{T}^2)})( 1+\|\varphi(s)\|_{H^{s,\gamma}_{g}}+\|\partial_{xy}^{s}U\|_{{L^2}(\mathbb{T}^2)})I,
\end{eqnarray}
where we have used Lemma \ref{y4.4} and Young's inequality.
Applying the classical maximum principle for parabolic equations, we have
\begin{eqnarray}
\begin{aligned}
&\|I(t)\|_{L^\infty(\mathbb{T}^2\times\mathbb{R}^+)}\\
&\leq \max \left\{            e^{C\left(1+G_1(t)\right)t}\|I(0)\|_{L^\infty(\mathbb{T}^2\times\mathbb{R}^+)},\mathop{\max}\limits_{\tau \in [0,t]}\left(e^{C\left(1+G_1(t)\right)(t-\tau)}\|I(\tau)\big{|}_{z=0}\|_{L^\infty(\mathbb{T}^2)}\right)
\right\}.
\label{3.36}
\end{aligned}
\end{eqnarray}
 To derive a lower bound estimate on $B_{(0,0,0)}$, we have
 \begin{eqnarray}
\begin{aligned}
&\left(\partial_t+u\partial_x+Ku\partial_{y}+(w+\frac{2\sigma}{1+z})\partial_z-\partial^2_z-\epsilon^2\partial^2_x-\epsilon^2\partial^2_y-\partial_{y}Ku\right)B_{(0,0,0)}\\
&=\left(\frac{\sigma}{1+z}w+\frac{\sigma(\sigma-1)}{(1+z)^2}\right)B_{(0,0,0)}.
\label{3.37}
\end{aligned}
\end{eqnarray}
Applying the classical maximum principle for (\ref{3.37}), we get
\begin{eqnarray}
\begin{aligned}
&\mathop{\min}\limits_{\mathbb{T}^2\times\mathbb{R}^+}(1+z)^{\sigma}\varphi(t)\\
&\geq \max \left\{         1-C\left(1+G_2(t)\right)te^{C\left(1+G_2(t)\right)t}
\right\}\times\min \left\{\mathop{\min}\limits_{\mathbb{T}^2\times\mathbb{R}^+}(1+z)^{\sigma}\varphi_0,        \mathop{\min}\limits_{[0,t]\times\mathbb{T}^2}\varphi\big{|}_{z=0}
\right\}.
\label{3.38}
\end{aligned}
\end{eqnarray}
Next we start estimating the boundary values, using  Lemma \ref{y4.3}, we obtain
\begin{eqnarray}
\begin{aligned}
\|I(\tau)|_{z=0}\|_{L^\infty(\mathbb{T})}&\leq C^2\sum \limits_{|\alpha|\leq 2}\big(\|D^{\alpha}\varphi\|_{L^2}^2+\|\partial_{z}D^{\alpha}\varphi\|_{L^2}^2+\|\partial_{xy}^{1}D^{\alpha}\varphi\|_{L^2}^2+\|\partial_{z}\partial_{xy}^{1}D^{\alpha}\varphi\|_{L^2}^2\\
&\quad +\|\partial_{xy}^{2}D^{\alpha}\varphi\|_{L^2}^2+\|\partial_{z}\partial_{xy}^{2}D^{\alpha}\varphi\|_{L^2}^2\big)\\
&\leq 6C^2\|\varphi\|_{H^{s,\gamma}_{g}}^{2},
\end{aligned}
\end{eqnarray}
which together with (\ref{3.36}) gives (\ref{33311}).

Now, we derive inequality (\ref{3.32}), according to (\ref{1.5}) and boundary condition $\partial^{\alpha_1}_x\partial^{\alpha_2}_y u=\partial^{\alpha_1}_x \partial^{\alpha_2}_y w=0$, we have
\begin{eqnarray}
\begin{aligned}
\partial_t D^{\alpha}\varphi\big{|}_{z=0}=D^{\alpha}(\epsilon^2\partial_x^2+\epsilon^2\partial_y^2+\partial_y^2 +\partial_{y}Ku-u\partial_x -Ku\partial_y-w\partial_z )\varphi\big{|}_{z=0}.
\end{aligned}
\end{eqnarray}
When $\alpha_3=0$ and $\alpha_1+\alpha_2 \leq 2$,
\begin{eqnarray}
\begin{aligned}
\partial_t D^{\alpha}\varphi\big{|}_{z=0}=(\epsilon^2\partial^{\alpha_1+2}_{x}\partial^{\alpha_2}_{y}+\epsilon^2\partial^{\alpha_1}_{x}\partial^{\alpha_2+2}_{y}
+\partial^{\alpha_1}_{x}\partial^{\alpha_2}_{y}\partial_z^2 )\varphi\big{|}_{z=0}.
\end{aligned}
\end{eqnarray}
When $1 \leq \alpha_3 \leq 2$ and $\alpha_1=\alpha_2=0$,
\begin{eqnarray}
\begin{aligned}
\partial_t D^{\alpha}\varphi\big{|}_{z=0}
&=\big(\epsilon^2\partial^{\alpha_3}_{z}\partial^{2}_{x}+\epsilon^2\partial^{\alpha_3}_{z}\partial^{2}_{y}+\partial^{\alpha_3+2}_{z}+(3\alpha_3-2)\partial_{y}K\varphi\partial_{z}^{\alpha-1}\\
&\quad-\alpha_3\varphi\partial_{z}^{\alpha_3-1}\partial_{x}-K\alpha_3\varphi\partial_{z}^{\alpha_3-1}\partial_{y}\big)\varphi\big{|}_{z=0}.
\end{aligned}
\end{eqnarray}
When $\alpha_3=1$ and $\alpha_1+\alpha_2 =1$,
\begin{eqnarray}
\begin{aligned}
\partial_t D^{\alpha}\varphi|_{z=0}&=(\epsilon^2\partial_{z}\partial_{xy}^{1}\partial^{2}_{x} +\epsilon^2\partial_{z}\partial_{xy}^{1}\partial^{2}_{y} +
\partial _{xy}^{1}\partial _{z}^{3}+\partial_{xy}^{1}\partial_{y}K\varphi+2\varphi\partial_{xy}+\varphi\partial_{xy}^{1}\partial_x\\
&\quad+\partial_{xy}^{1}\varphi\partial_{x}+K\varphi\partial_{xy}^{1}\partial_y+\partial_{xy}^{1}(K\varphi)\partial_{y}
)\varphi\big{|}_{z=0}.
\end{aligned}
\end{eqnarray}
For $s\geq 7$, using Lemma \ref{y4.4}, we have
 \begin{eqnarray}
\begin{aligned}
\|\partial_t I\big{|}_{z=0}\|_{{L^\infty}(\mathbb{T}^2)}&\leq C_{s, \gamma  } \|D^{\alpha}\varphi\partial_t D^{\alpha}\varphi \big{|}_{z=0}\|_{{L^\infty}(\mathbb{T}^2)}\\
&\leq  C_{s, \gamma  }\|\varphi(\tau)\|_{H^{s,\gamma}_{g}}(\|\varphi(\tau)\|_{H^{s,\gamma}_{g}}+\|\partial_{xy}K(\tau)\|_{{L^{\infty}}(\mathbb{T}^2)}\|\varphi(\tau)\|_{H^{s,\gamma}_{g}}^2).
\end{aligned}
\end{eqnarray}
Hence, a direct integration yields
 \begin{eqnarray}
\begin{aligned}
\|I(t)\big{|}_{z=0}\|_{{L^\infty}(\mathbb{T}^2)}\leq  \|I(0)\big{|}_{z=0}\|_{{L^\infty}(\mathbb{T}^2)}+C_{s, \gamma  } [1+Y(t)W(t)]W(t)^2t,
\end{aligned}
\end{eqnarray}
which, together with (\ref{3.36}), gives (\ref{3.32}). For $s\geq4$,
 \begin{eqnarray}
\begin{aligned}
\|\partial_t \varphi|_{z=0}\|_{{L^\infty}(\mathbb{T}^2)}\leq  CW(t).
\end{aligned}
\end{eqnarray}
Thus, a direct integration leads to
 \begin{eqnarray}
\begin{aligned}
\mathop{\min}\limits_{\mathbb{T}^2}{\varphi(t)}|_{z=0}\geq \mathop{\min}\limits_{\mathbb{T}^2}{\varphi_0}\big{|}_{z=0}-CW(t)t,
\end{aligned}
\end{eqnarray}
which, together with (\ref{3.38}), gives (\ref{3.33}). The proof of Lemma \ref{y2.3} is thus complete.
$\hfill$$\Box$
\section{Local-in-time existence and uniqueness}
\subsection{Local-in-time existence }
In this subsection, we go back to use the symbol $(u^\epsilon, K^\epsilon u^{\epsilon}, w^\epsilon, \varphi^{\epsilon})$ instead of $(u, Ku, w, \varphi)$ from Subsections 2.1-2.4 to denote the solution to the regularized system (\ref{1.5}). To obtain the local-in-time solution of the initial-boundary value problem (\ref{1.1})-(\ref{1.2}) or (\ref{1.4})-(\ref{1.44}), we will construct the solution to the Prandtl equations (\ref{1.4}) by passing to the limit $\epsilon \rightarrow 0^+$ in the regularized Prandtl equations (\ref{1.5}). We only sketch the proof into five steps and more details can be found in \cite{3}.\\
$\textbf{Step 1.}$ According to the definition of $F$, assumption (\ref{1.8}) and (\ref{1.88}), and the regularized Bernoulli's law,
\begin{eqnarray}
\|F\|_{L^\infty} \leq C_{s, \gamma,\sigma,\delta  }\left(M_{U}+M_{K}\right) < +\infty,
\end{eqnarray}
which, together with (\ref{3.2222}), implies
\begin{eqnarray}
\begin{aligned}
\|\varphi\|_{H^{s,\gamma}_{g}}^2 &\leq \left\{\|\varphi_{0}\|_{H^{s,\gamma}_{g}}^2+C_{s, \gamma,\sigma,\delta  }\left(M_{U}+M_{K}\right)t_{1}\right\}\\
&\quad \times \left\{1-C_{s, \gamma,\sigma,\delta  }\left(\|\varphi_{0}\|_{H^{s,\gamma}_{g}}^2+C_{s, \gamma,\sigma,\delta  }\left(M_{U}+M_{K}\right)t_{1}\right)^{s-1}t_{2}\right\}^{-\frac{1}{s-1}},
\label{2.01555}
\end{aligned}
\end{eqnarray}
for $t_{1}$ and $t_{2}$ to be chosen later. If choosing $t_{1} \leq \frac{3\|\varphi_0\|_{H^{s,\gamma}_{g}}^2}{C_{s, \gamma,\sigma,\delta  }(M_{U}+M_{K})}$ in (\ref{2.01555}), then we can  get
\begin{eqnarray}
\begin{aligned}
\|\varphi\|_{H^{s,\gamma}_{g}}^2 &\leq 4\|\varphi_{0}\|_{H^{s,\gamma}_{g}}^{2}\times \left\{1-C_{s, \gamma,\sigma,\delta  }\left(4\|\varphi_{0}\|_{H^{s,\gamma}_{g}}^{2}\right)^{s-1}t_{2}\right\}^{-\frac{1}{s-1}}.
\label{2.015555}
\end{aligned}
\end{eqnarray}
Further, choosing $t_{2} \leq \frac{1-2^{2-2s}}{2^{2s-2}C_{s, \gamma,\sigma,\delta  }\|\varphi_0\|_{H^{s,\gamma}_{g}}^{2s-2}}$ in (\ref{2.015555}) gives
\begin{eqnarray}
\begin{aligned}
\|\varphi\|_{H^{s,\gamma}_{g}}^2 &\leq 4\|\varphi_{0}\|_{H^{s,\gamma}_{g}}^{2}\times \left(2^{2-2s}\right)^{-\frac{1}{s-1}} \leq 16\|\varphi_{0}\|_{H^{s,\gamma}_{g}}^{2}.
\end{aligned}
\end{eqnarray}
Note that the functions in both brackets in  (\ref{2.01555}) are monotonically increasing with respect to time $t_{1}$ and $t_{2}$, thus we derive the uniform estimate for any $\epsilon \in [0,1]$ and any $t \in [0,T_1]$,
\begin{eqnarray}
\|\varphi^\epsilon\|_{H^{s,\gamma}_{g}} \leq 4 \|\varphi_0\|_{H^{s,\gamma}_{g}},
\label{2.0155555}
\end{eqnarray}
provided that $T_1$ is chosen  such that
\begin{eqnarray*}
T_1 :=\min \left\{\frac{3\|\varphi_0\|_{H^{s,\gamma}_{g}}^2}{C_{s, \gamma,\sigma,\delta  }(M_{U}+M_{K})},\frac{1-2^{2-2s}}{2^{2s-2}C_{s, \gamma,\sigma,\delta  }\|\varphi_0\|_{H^{s,\gamma}_{g}}^{2s-2}}\right\}.
\end{eqnarray*}
$\textbf{Step 2.}$ When $s \geq 7$, we know from definition (\ref{2.0122}) of $W$ and (\ref{2.012}) of $G_1$ that for any $t \in [0,T_1]$,
\begin{eqnarray}
W(t) \leq 4 \|\varphi_0\|_{H^{s,\gamma}_{g}} \quad and \quad G_1(t) \leq M_{K}+(4 \|\varphi_0\|_{H^{s,\gamma}_{g}}+M_{U})(1+M_{K}).
\label{2.015}
\end{eqnarray}
Thus, similar to the selection method of $T_{1}$, if we choose
\begin{eqnarray*}
T_2 :=\min \left\{T_1,\frac{1}{64\delta^2 C_{s, \gamma }(1+4 \|\varphi_0\|_{H^{s,\gamma}_{g}}) \|\varphi_0\|_{H^{s,\gamma}_{g}}^2},\frac{\ln 2}{C_{s, \gamma,\sigma,\delta  }[1+M_{K}+(4 \|\varphi_0\|_{H^{s,\gamma}_{g}}+M_{U})(1+M_{K})]}\right\},
\end{eqnarray*}
then using inequality (\ref{3.32}) and initial assumption
\begin{eqnarray}
\sum \limits_{|\alpha|\leq 2}  |(1+z)^{\sigma+\alpha_{3}}D^{\alpha}\varphi_0 | ^{2}\leq\frac{1}{4\delta^{2}},
\end{eqnarray}
we have the upper bound
\begin{eqnarray}
\Big\|\sum \limits_{|\alpha|\leq 2}  |(1+z)^{\sigma+\alpha_{3}}D^{\alpha}\varphi^\epsilon (t) | ^{2}\Big\|_{L^{\infty}(\mathbb{T}^2\times \mathbb{R}^+)}\leq\frac{1}{\delta^{2}}
\label{2.019}
\end{eqnarray}
for all $t \in [0,T_2]$. When $s \geq 5$, from the initial hypothesis $\|\varphi_0\| \leq C\delta^{-1}$, we have the same estimate (\ref{2.019}) for all $t\in [0,T_2]$.\\
$\textbf{Step 3.}$
Let us choose
\begin{eqnarray*}
T_3 :=\min \left\{T_1,\frac{\delta}{8 C_{s, \gamma }\|\varphi_0\|_{H^{s,\gamma}_{g}}},\frac{1}{6C_{s, \gamma,\sigma,\delta  }N},\frac{\ln 2}{C_{s, \gamma,\sigma,\delta  }N}\right\},
\end{eqnarray*}
where $N :=1+M_{K}+(4 \|\varphi_0\|_{H^{s,\gamma}_{g}}+M_{U})(1+M_{K})$. Then by using (\ref{2.012}) and (\ref{2.015}), we derive the uniform estimate for any $\epsilon \in [0,1]$ and any $t \in [0,T_3]$,
\begin{eqnarray}
\min \limits_{\mathbb{T}^2\times \mathbb{R}^+}(1+z)^{\sigma} \varphi^\epsilon(t) \geq \delta.
\end{eqnarray}
$\textbf{Step 4.}$
In summary, the above uniform estimates  hold for any $t \in [0,T]$, if $T$ is chosen to satisfy $T=\min\{T_1,T_2,T_3\}$. Further using almost equivalence relation  (\ref{4.1}), Lemma \ref{y4.4} and (\ref{2.0155555}), we have
\begin{align}
\sup  \limits_{0 \leq t \leq T}( \|\varphi^\epsilon\|_{H^{s,\gamma}}+\|u^\epsilon-U\|_{H^{s,\gamma-1}}) &\leq 4C\|\varphi_0\|_{H^{s,\gamma}_{g}}+ C\sup \limits_{0 \leq t \leq T}\left(\sum_{k=0}^{s}\|(1+z)^{\gamma-1}\partial_{xy}^{k}(u^\epsilon-U)\|_{L^{2}}\right)\nonumber\\
&\leq C\|\varphi_0\|_{H^{s,\gamma}_{g}}+ C_{s, \gamma,\sigma,\delta  }\sup \limits_{0 \leq t \leq T}\left(\|\varphi\|_{H^{s,\gamma}_{g}}+\|\partial_{xy}^s U\|_{L^2(\mathbb{T}^2)}\right)\nonumber\\
& \leq C_{s, \gamma,\sigma,\delta  }\left(\|\varphi_0\|_{H^{s,\gamma}_{g}}+\sup \limits_{0 \leq t \leq T}\|\partial_{xy}^s U\|_{L^2(\mathbb{T}^2)}\right)< + \infty.
\label{2.016}
\end{align}
From the equation (\ref{1.4}),  (\ref{1.5}), (\ref{2.016}) and Lemma \ref{y4.4}, we also have $\partial_t \varphi^{\epsilon}$	 and $\partial_t (u^{\epsilon}-U)$ are uniformly
bounded in $L^{\infty}([0,T];H^{s-2,\gamma})$ and $L^{\infty}([0,T];H^{s-2,\gamma-1})$ respectively. By Lions-Aubin lemma and the compact embedding of $H^{s,\gamma}$ in $H^{s^{\prime}}_{loc}$
, we have after taking a subsequence, as $\epsilon_k \rightarrow 0^+$,
\begin{equation}\left\{
\begin{array}{ll}
\varphi^{\epsilon_k} \stackrel{*}{\rightharpoonup} \varphi, &{\rm in} \quad L^{\infty}([0,T];H^{s,\gamma}),\\
\varphi^{\epsilon_k} \rightarrow \varphi, &{\rm in} \quad C([0,T];H^{s^{\prime}}_{loc}),\\
u^{\epsilon_k}-U \stackrel{*}{\rightharpoonup} u-U,  &{\rm in} \quad L^{\infty}([0,T];H^{s,\gamma-1}),\\
u^{\epsilon_k} \rightarrow u, &{\rm in} \quad C([0,T];H^{s^{\prime}}_{loc}),
\end{array}
         \right.\end{equation}
for all $s^{\prime} < s$, where
\begin{equation}\left\{
\begin{array}{ll}
\varphi=\partial_z u \in L^{\infty} ([0,T];H^{s,\gamma})\cap \bigcap_{s^\prime <s}C([0,T];H^{s^{\prime}}_{loc}),\\
 u-U \in L^{\infty} ([0,T];H^{s,\gamma-1})\cap \bigcap_{s^\prime <s}C([0,T];H^{s^{\prime}}_{loc}).
\end{array}
         \right.\end{equation}
$\textbf{Step 5.}$
Using the local uniform
convergence of $\partial_x u^{\epsilon_k}$ and $\partial_y (K^{\epsilon_k}u ^{\epsilon_k})$, we also have the pointwise convergence of $w^{\epsilon_k}:$ as $\epsilon \rightarrow 0^+$,
\begin{eqnarray}
w^{\epsilon_k}=-\int^{z}_{0} \partial_x u ^{\epsilon_k}~ dz-\int^{z}_{0} \partial_y (K^{\epsilon_k}u ^{\epsilon_k})~ dz \rightarrow -\int^{z}_{0} \partial_x u ~ dz-\int^{z}_{0} \partial_y (Ku )~ dz=: w.
\end{eqnarray}
Thus, passing the limit $\epsilon_k \rightarrow 0^+$ in the initial-boundary value problem (\ref{1.4}), we get that the limit
$(u, Ku, w)$ solves the initial-boundary value problem (\ref{1.1}) in the classical sense. This completes the proof
of the existence.
\subsection{Uniqueness of solutions}
In this subsection, we are going to prove the uniqueness of solutions to 3D  Prandtl model.  Let us denote $(\bar u, K\bar{u}, \bar{w})=(u_1, Ku_1, w_1)-(u_2, Ku_2, w_2)$, $\bar {\varphi} ={\varphi}_{1}-{\varphi}_{2}$, $a_2=\frac{\partial_z \varphi_2}{\varphi_{2}}$ and $\bar g =\bar {\varphi} -a_2 \bar u$. It is easy to check that $\bar g =\bar {\varphi} -a_2 \bar u=\varphi_2 \partial_z(\frac{\bar u}{\varphi_2})$ and the evolution equation on $\bar g$ is as follows
\begin{eqnarray}
\begin{aligned}
&(\partial_t+u_1\partial_x+Ku_1\partial_y+w_1\partial_z-\partial^{2}_{z})\bar g\\
&=(\partial_t+u_1\partial_x+Ku_1\partial_y+w_1\partial_z-\partial^2_z)\bar \varphi-a_2(\partial_t+u_1\partial_x+Ku_1\partial_y+w_1\partial_z-\partial^2_z)\bar u\\
&\quad -\bar u(\partial_t+u_1\partial_x+Ku_1\partial_y+w_1\partial_z-\partial^2_z)a_2-2\partial_z a_2 \bar {\varphi}.
 \label{33.40}
\end{aligned}
\end{eqnarray}
Next, we calculate the values of the first three terms on the right-hand  side of the equality (\ref{33.40}) respectively, recalling the vorticity system (\ref{1.5}), we have
\begin{eqnarray}
\begin{aligned}
(\partial_t+u_i\partial_x+Ku_i\partial_y+w_i\partial_z)\partial_z \varphi_i&=\partial_z^3 \varphi_i+[\partial_x u_i+\partial_y (Ku_i)] \partial_{z} \varphi_{i} -\varphi_i \partial_x \varphi_i-K\varphi_i \partial_y \varphi_i\\
&\quad+\partial_{y}K(\partial_{z} \varphi u+\varphi^{2}),
\end{aligned}
\end{eqnarray}
then, according to the definition of $a_i$, we get
\begin{align}
(\partial_t+u_i\partial_x+Ku_i\partial_y+w_i\partial_z)a_i&=\frac{(\partial_t+u_i\partial_x+Ku_i\partial_y+w_i\partial_z)\partial_z \varphi_i}{\varphi_i}-\frac{\partial_z \varphi_i(\partial_t+u_i\partial_x+Ku_i\partial_y+w_i\partial_z) \varphi_i}{\varphi_i^2}\nonumber\\
&=\frac{\partial^3_z \varphi_i}{\varphi_i}+a_i\partial_x u_i+a_i\partial_y (Ku_i)-\partial_x \varphi_i-K\partial_y \varphi_i -a_i \frac{\partial_z^2 \varphi_i}{\varphi_i}+\partial_{y}K\varphi_{i},
\end{align}
which, combined with the fact
\begin{eqnarray*}
\frac{\partial^3_z \varphi_i}{\varphi_i}-a_i \frac{\partial_y^2 a_i}{\varphi_i}=\partial_z^2 a_i+2a_i \partial_z a_i,
\end{eqnarray*}
implies
\begin{eqnarray}
\begin{aligned}
(\partial_t+u_i\partial_x+Ku_i\partial_y+w_i\partial_z-\partial_{z}^{2})a_i=a_i\partial_x u_i+a_i\partial_y (Ku_i)-\partial_x \varphi_i-K\partial_y \varphi_i +\partial_{y}K\varphi +2a_i \partial_z a_i.
\end{aligned}
\end{eqnarray}
Furthermore, we conclude that
\begin{eqnarray}
\begin{aligned}
&(\partial_t+u_1\partial_x+Ku_1\partial_y+w_1\partial_z-\partial_{z}^{2})a_2\\
&=a_2\partial_x u_2+a_2\partial_y (Ku_2)-\partial_x \varphi_2-K\partial_y \varphi_2 +\partial_{y}K\varphi_{2} +2a_2 \partial_z a_2+(\bar u \partial_x+K\bar {u} \partial_y+\bar{w}\partial_{z})a_2.
 \label{33.41}
\end{aligned}
\end{eqnarray}
For the estimate of $(\partial_t+u_1\partial_x+Ku_1\partial_y+w_1\partial_z-\partial_{z}^{2})\bar u$ and $(\partial_t+u_1\partial_x+Ku_1\partial_y+w_1\partial_z-\partial_{z}^{2})\bar \varphi$, from the evolution equation on $u$ and on $\varphi$, we can derive that
\begin{eqnarray}
\begin{aligned}
(\partial_t+u_1\partial_x+Ku_1\partial_y+w_1\partial_z-\partial_{z}^{2})\bar u&=-\bar u \partial_x u_2-K\bar {u} \partial_y  u_2-\bar {w} \partial_z  u_2,
\end{aligned}
\end{eqnarray}
and
\begin{eqnarray}
\begin{aligned}
(\partial_t+u_1\partial_x+Ku_1\partial_y+w_1\partial_z-\partial_{z}^{2})\bar {\varphi}&=-\bar u \partial_x \varphi_2-K\bar {u} \partial_y  \varphi_2-\bar {w} \partial_z  \varphi_2+\partial_{y}K(\bar{u}\varphi_{2}+u_1\bar{\varphi}).
 \label{33.42}
\end{aligned}
\end{eqnarray}
Using (\ref{33.40}) and (\ref{33.41})-(\ref{33.42}), we have
\begin{eqnarray}
\begin{aligned}
(\partial_t+u_1\partial_x+Ku_1\partial_y+w_1\partial_z-\partial_{z}^{2})\bar g&= -\bar u(\bar u \partial_xa_{2}+K\bar {u} \partial_ya_{2}+\bar{w}\partial_{z}a_{2}+2a_{2}\partial_{z}a_{2})-2\bar {\varphi} \partial_z a_2\\
&\quad+\partial_y{K}(u_1\bar{\varphi}-\bar{u}a_{2}u_{2}).
\end{aligned}
\end{eqnarray}
Now we derive $L^2$ estimate on $\bar g$. For any $t \in (0,T]$, multiplying by $2 \bar g$ and then integrating by parts over $\mathbb{T}^2\times\mathbb{R}_{+}$, we obtain
\begin{eqnarray}
\begin{aligned}
&\frac{d}{dt}\|\bar g\|_{L^2}^2+2\|\partial_z \bar g\|_{L^2}^2\\
&=\iint_{\mathbb{T}^{2}}\bar g \partial_z \bar g \big{|}_{z=0}dxdy-2\iiint \bar g \bar u(\bar u \partial_xa_{2}+K\bar {u} \partial_ya_{2}+\bar{w}\partial_{z}a_{2}+2a_{2}\partial_{z}a_{2})-4\iiint \bar g \bar {\varphi} \partial_z a_2  \\
&\quad+2\iiint \bar{g}\partial_y{K}(u_1\bar{\varphi}-\bar{u}a_{2}u_{2}) -2\iiint \bar g(u_1\partial_x\bar{g}+Ku_1\partial_y\bar{g}+w_1\partial_z\bar{g}).
 \label{33.43}
\end{aligned}
\end{eqnarray}
We need to estimate all terms on the right-hand side of equation (\ref{33.43}). Applying the simple trace theorem and Young's inequality,
\begin{align}
\left|\int_{\mathbb{T}^2}\bar {g} \partial_z \bar {g} |_{z=0}dxdy\right|
&\leq|\iint a_2 |\bar g|^2  dxdydz|+|\iint \partial_z a_2 |\bar g|^2  dxdydz|+2|\iint  a_2 \bar g \partial_z \bar g  dxdydz|\nonumber\\
&\leq \frac{1}{2}\|\partial_z \bar g\|_{L^2}^2+C_{\sigma,\delta  }\|\bar g\|_{L^2}^2,
\label{w3.20}
\end{align}
where we have used the fact $ \partial_z \bar g \big{|}_{z=0}=-a_2 \bar g  |_{z=0}$ . We claim $\|\frac{\bar u}{1+z}\|_{L^2} \leq C_{\sigma,\delta} \|\bar g\|_{L^2}$, so by Lemma \ref{y4.4},
\begin{align}
&-2\iint \bar g \bar u(\bar u \partial_x a_2+\bar v \partial_y a_2+2a \partial_y a_2) \nonumber\\
&\leq  C\|(1+z)(\bar u \partial_xa_{2}+K\bar {u} \partial_ya_{2}+\bar{w}\partial_{z}a_{2}+2a_{2}\partial_{z}a_{2})\|_{L^\infty}\left\|\frac{\bar u}{1+z}\right\|_{L^2}\|\bar g \|_{L^2}\nonumber\\
&\leq C_{ \gamma,\sigma,\delta  }\left(1+\|w_i\|_{H^{5,\gamma}_{g}}+\|\partial_{xy}^{5}U\|_{{L^2}(\mathbb{T}^2)}+\|K\|_{{L^{\infty}}(\mathbb{T}^2)}\|w_i\|_{H^{5,\gamma}_{g}}
+\|K\|_{{L^{\infty}}(\mathbb{T}^2)}\|\partial_{xy}^{5}U\|_{{L^2}(\mathbb{T}^2)}\right)\|\bar g\|_{L^2}^2.
\end{align}
Now we shall prove  that $\|\frac{\bar u}{1+z}\|_{L^2}$ can be controlled by $\|\bar g\|_{L^2}$. Indeed, since $\delta \leq (1+z)^{\delta} \varphi_2\leq \delta^{-1}$ and Lemma \ref{y4.1},
\begin{eqnarray}
\begin{aligned}
\left\|\frac{\bar u}{1+z}\right\|_{L^2} \leq \delta^{-1}\left \|(1+z)^{-\sigma-1}\frac{\bar u}{\varphi_2}\right\|_{L^2} \leq C_{\sigma,\delta}\left\|(1+z)^{-\sigma} \partial_{z}\left(\frac{\bar u}{\varphi_2}\right)\right\|_{L_2} \leq C_{\sigma,\delta} \|\bar g\|_{L^2}.
\end{aligned}
\end{eqnarray}
In addition, we also have
\begin{eqnarray}
\begin{aligned}
\| \bar \varphi\|_{L^2} \leq \|\bar g \|_{L^2}+\delta^{-2}\left\|\frac{\bar u}{1+z}\right\|_{L^2}\leq C_{\sigma,\delta} \|\bar g\|_{L^2},
\end{aligned}
\end{eqnarray}
then we can obtain that
\begin{eqnarray}
\begin{aligned}
-4\iint  \bar g \bar \varphi \partial_z a_2 \leq C_{\sigma,\delta} \|\bar g\|_{L^2}^2
\end{aligned}
\end{eqnarray}
and
\begin{eqnarray}
\begin{aligned}
2\iiint \bar{g}\partial_y{K}(u_1\bar{\varphi}-\bar{u}a_{2}u_{2})\leq C_{ \gamma,\sigma,\delta  }\|\partial_{y}K\|_{{L^{\infty}}(\mathbb{T}^2)}\left(\|w_i\|_{H^{5,\gamma}_{g}}
+\|\partial_{xy}^{5}U\|_{{L^2}(\mathbb{T}^2)}\right)\|\bar g\|_{L^2}^2.
\end{aligned}
\end{eqnarray}
For the last term in (\ref{33.43}), using the integration by parts, boundary condition $(u_1, w_1)|_{z=0}=0$ and $\partial_x u_1+\partial_y (Ku_1)+\partial_z w_1=0$, we readily  show
\begin{eqnarray}
\begin{aligned}
-2\iiint \bar g(u_1\partial_x\bar{g}+Ku_1\partial_y\bar{g}+w_1\partial_z\bar{g})=0.
\label{w3.26}
\end{aligned}
\end{eqnarray}
Combining all of the above estimates (\ref{w3.20})-(\ref{w3.26}), we can derive from (\ref{33.43})
\begin{eqnarray}
\begin{aligned}
\frac{d}{dt}\|\bar g\|_{L^2}^2 \leq C_{ \gamma,\sigma,\delta  }\left[1+\left(1+\|\partial_{y}K\|_{{L^{\infty}}(\mathbb{T}^2)}\right)\left(\|w_i\|_{H^{5,\gamma}_{g}}+\|\partial_{xy}^{5}U\|_{{L^2}(\mathbb{T}^2)}\right)\right]\|\bar g\|_{L^2}^2,
\end{aligned}
\end{eqnarray}
which, together with Gronwall's inequality, yields
\begin{eqnarray}
\begin{aligned}
\|\bar g(t)\|_{L^2}^2 \leq \|\bar g(0)\|_{L^2}^2e^{Ct},
\end{aligned}
\end{eqnarray}
here  $C=C_{ \gamma,\sigma,\delta  }\left[1+\left(1+\|\partial_{y}K\|_{{L^{\infty}}(\mathbb{T}^2)}\right)\left(\|w_i\|_{H^{5,\gamma}_{g}}+\|\partial_{xy}^{5}U\|_{{L^2}(\mathbb{T}^2)}\right)\right]$, and this implies $\bar g=0$ due to $u_1|_{t=0}=u_2|_{t=0}$. Since $\varphi_{2}\partial_z (\frac{u_1-u_2}{\varphi_{2}})=\bar g =0$, we have
\begin{eqnarray}
\begin{aligned}
u_1-u_2=q \varphi_{2}
\end{aligned}
\end{eqnarray}
for some function $q=q(t,x)$. By using the  Oleinik's monotonicity assumption $\varphi_{2} >0 $ and boundary condition $u_1|_{z=0}=u_2|_{z=0}=0$, we can get $q=0$, and hence $u_1=u_2$. Furthermore, using $\partial_{x} u+\partial_{y}(K u)+\partial_{z} w=0$, then $w_i$   can be uniquely determined (i.e., $w_1=w_2$ ). This completes the proof of the uniqueness of solutions.

$\mathbf{Acknowledgment}$

This paper was supported in part by the NNSF of China with contract  number 12171082, the fundamental research funds for the central universities with contract numbers 2232022G-13, 2232023G-13 and  the Scientific Research Foundation of Education Department of Yunnan province with  contract  number 2023J0133.

$\mathbf{Conflict\,\,of\,\,interest\,\,statement}$

The authors have no conflict of interest.

\end{document}